\documentclass[a4paper,times,3p]{elsarticle}
\usepackage[utf8]{inputenc}

\usepackage[T1]{fontenc}
\usepackage{microtype}

\usepackage{booktabs}

\usepackage{xcolor}
\usepackage{amsmath}
\usepackage{amssymb}
\usepackage{amsthm}
\usepackage{bm}
\usepackage{graphicx}
\usepackage{subfigure}
\usepackage{epstopdf}
\usepackage{psfrag}
\usepackage{color}
\usepackage{hyperref}
\usepackage{transparent}

\usepackage{mathtools}

\definecolor{cmykcyan}{cmyk}{1,0,0,0}
\definecolor{cmykred}{cmyk}{0,1,1,0}
\definecolor{cmykblack}{cmyk}{0,0,0,1}

\newtheorem*{remark}{Remark}

\usepackage[font=footnotesize,labelfont=bf]{caption}
\usepackage{multirow}

\usepackage{pxfonts}

\usepackage{import}

\graphicspath{./img/}

\usepackage[textsize=scriptsize]{todonotes}
\usepackage{setspace}

\journal{International Journal for Numerical Methods in Fluids}

\begin{document}

\begin{frontmatter}

\title{Combining Boundary-Conforming Finite Element Meshes on Moving Domains Using a Sliding Mesh Approach}

\author[add1]{Jan Helmig\corref{cor1}}
\ead{helmig@cats.rwth-aachen.de}
\author[add1,add2]{Fabian Key\corref{}}
\ead{key@cats.rwth-aachen.de}
\author[add1]{Marek Behr\corref{}}
\ead{behr@cats.rwth-aachen.de}
\author[add1,add2]{Stefanie Elgeti\corref{}}
\ead{elgeti@ilsb.tuwien.ac.at}

%\ead[url]{http://www.cats.rwth-aachen.de}
\cortext[cor1]{Corresponding author}

\address[add1]{Chair for Computational Analysis of Technical Systems (CATS)\\
	RWTH Aachen University, 52056 Aachen, Germany}
\address[add2]{Institute of Lightweight Design and Structural Biomechanics\\
	TU Wien, Getreidemark 9 A-1060 Vienna, Austria}

\begin{abstract}
For most finite element simulations, boundary-conforming meshes have significant advantages in terms of accuracy or efficiency.
This is particularly true for complex domains. However, with increased complexity of the domain, generating a boundary-conforming mesh becomes more difficult and time consuming.
One might therefore decide to resort to an approach where individual boundary-conforming meshes are pieced together in a modular fashion to form a larger domain.
This paper presents a stabilized finite element formulation for fluid and temperature equations on sliding meshes.
It couples the solution fields of multiple subdomains whose boundaries slide along each other on common interfaces. Thus, the method allows to use highly tuned boundary-conforming meshes for each subdomain that are only coupled at the overlapping boundary interfaces. In contrast to standard overlapping or fictitious domain methods the coupling is broken down to few interfaces with reduced geometric dimension.
The formulation consists of the following key ingredients: the coupling of the solution fields on the overlapping surfaces is imposed weakly using a stabilized version of  Nitsche's method. It ensures mass and energy conservation at the common interfaces. Additionally, we allow to impose weak Dirichlet boundary conditions at the non-overlapping parts of the interfaces.
We present a detailed numerical study for the resulting stabilized formulation. It shows optimal convergence behavior for both Newtonian and generalized Newtonian material models. Simulations of flow of plastic melt inside single-screw as well as twin-screw extruders demonstrate the applicability of the method to complex and relevant industrial applications.
\end{abstract}

\end{frontmatter}

% !TEX root = ./main.tex
%!TeX spellcheck = en-US
\section{Introduction}

In this paper we present a stabilized finite element formulation for fluid and temperature equations that allows to couple boundary-conforming discretizations of individual moving domains at common interfaces. In the world of finite element analysis for flow problems, the representation of the computational domain plays an important role, especially for complex moving domains in 3D. The numerical solution requires high-quality meshes to ensure good approximation properties on the one hand, and a proper geometric resolution of the given domain on the other hand. The question of how to balance these two aims has led to a variety of approaches and methods. Broadly speaking, these can be placed into one of two general categories: (1) boundary-conforming meshes that are aligned with the domain boundary and (2) unfitted, fictitious or overlapping methods, which all describe techniques in which the actual domain is embedded into a static background mesh or individual meshes arbitrarily overlap. \\
In the former case, a mesh has to be generated based on given geometrical data. This can be an expensive task both in terms of human resource –- since manual intervention is often necessary –- as well as in terms of computer resources. At the same time, this approach ensures full control over the resolution and can easily include expert knowledge, e.g., it is possible to generate high-quality boundary layers or to refine special regions of interest already in advance. Furthermore, the imposition of Dirichlet boundary conditions is straightforward. In case of moving domains, e.g., for fluid-structure interaction the mesh can be updated in order to adapt to the moving boundary. What may seem like a drawback is, at the current state of the art, easily covered using mesh update methods based on radial basis functions \cite{de2007mesh,rendall2009efficient} or the Elastic Mesh Update Method (EMUM) \cite{johnson1994mesh,stein2003mesh}.
A more advanced extension of EMUM based on fiber-reinforced hyperelasticity is pesented in \cite{takizawa2020low}. In addition, there exist a broad range of specialized mesh update methods for specific applications such as the Shear Slip Mesh Update Method (SSMUM) for rotating components \cite{behr1999shear,behr2003shear} or its extension to large translation, the Virtual Ring Shear Slip Update Method (VRSSMUM) \cite{key2018virtual}.
For rotating screw machines, methods that automatically adapt to the moving screw domain have been developed in \cite{helmig2019boundary,hinz2019boundary,rane2013grid}.
Yet, even the most sophisticated method has its limitations: They occur when  boundary deformations are too large or result in topological changes. At this point, the only standard option that is available is continuous remeshing during the simulation; the price to pay is the extensive use of computational resources along with a loss of accuracy of the simulation due to the necessary mapping between the individual meshes. \\
With regard to unfitted methods, a variety of options such as level-set methods \cite{lehrenfeld2017higher} or immersed boundary methods \cite{peskin2002immersed} are available. In the following, we will focus on methods using overlapping meshes, as these are relevant to the further development of the paper. For overlapping mesh methods, meshes are generated around static or moving objects inside the domain –- e.g., the rotating object –- and are then embedded into a static background mesh that covers the whole computational domain. They have been introduced under the name of Chimera \cite{steger1991chimera} or Overset methods \cite{belk1995role}. The solutions on the domains are coupled weakly using a Dirichlet/Neumann coupling \cite{houzeaux2003chimera}.
The key advantage of simple mesh generation is somewhat offset by the disadvantage of the resulting in additional coupling steps. A general overview of these methods can be found in \cite{houzeaux2017domain}.
Fictitious domain methods using Lagrange multipliers that enrich the finite element function space (XFEM) are also wide spread, and have been applied to numerous complex applications \cite{atamian1991control,gerstenberger2008extended,mayer20103d,fard2012extended,villanueva2014density}.
Instead of enriching the functions space, Nitsche’s method \cite{nitsche1971variationsprinzip} has been used in \cite{hansbo2002unfitted} to enforce the interface coupling as well as boundary conditions weakly. These Nitsche-based methods form the basis of cutFEM, since two meshes intersect or are cut by each other, see, e.g., in \cite{massing2014stabilized,schott2016stabilized}.
Based on this, a multimesh approach has been developed that allows to couple arbitrary many overlapping meshes \cite{johansson2019multimesh,dokken2019multimesh}. Noteworthy are also formulations of the finite cell method that in addition use Nitsche’s method \cite{parvizian2007finite,hoang2017mixed}. In terms of rotating screw machines, standard fictitious domain methods accompanied with a mesh deformation technique have been developed in \cite{ianus2014mesh}.
All the fictitious domain methods have in common that they require the computation of the cut between either the background mesh with the underlying geometry or the overlapping meshes. This is not trivial in 3D. Furthermore, load balancing is challenging for highly parallel large scale computations. Bazilevs and Hughes present an approach that combines the best from both worlds for cases including rotating with prescribed rotation \cite{bazilevs2008nurbs}. An individual, boundary-conforming mesh is generated both around the rotating object and for the outer domain. In contrast to cutFEM, both domains do not overlap but share a common interface. The two meshes slide over this interface throughout the rotation, resulting in a non-matching interface one dimension lower than the actual domain. The two solution fields are coupled using Nitsche’s method similar to cutFEM. The method has been successfully applied to large-scale simulations of wind turbines \cite{hsu2012fluid}.
Note, however, that this method does not handle interfaces that change over time. \\
The method that will be presented in this work follows Bazilevs and Hughes \cite{bazilevs2008nurbs} in the sense that it allows for boundary-conforming meshes wherever possible and couples those meshes at common interfaces. As before, these interfaces may change throughout the simulation when the meshes slide along the interface. Like Bazilevs and Hughes, we only consider overlapping meshes at boundary interfaces –  e.g., a 2D cross section in 3D – thus strongly simplifying the computation of cut cells and avoiding degenerated polyhedrons. What is new is that we also take into account non-matching interfaces caused by rigid structures that exist only one side of the common interface. Here, we make use of weak imposition of boundary conditions in order to match the idea of the coupling condition \cite{bazilevs2007weak,massing2018stabilized}. Furthermore, we have expanded the formulation to generalized Newtonian material models and to conjugate heat transfer. The method is especially suited for applications, where individual complex moving domains are modularly stacked together in various orders. In this paper, we consider the example of extruders. \\
The paper is structured as follows: In Section 2, we give a more detailed motivation for the presented method as well as the problem setting of the sliding mesh approach. Moreover, the governing equations as well as the coupling conditions are formulated. Section 3 states the weak stabilized formulation for each domain. Additionally, a detailed description of the weak imposition of the coupling conditions and boundary conditions using Nitsche’s method is presented. This is followed by a discussion of the practical implementation. Numerical examples are presented in Section 4. Spatial convergence analysis for viscous and convective flow regimes for the 2D Taylor-Green flow, as well as for a complex 3D flow of plastic melt through a twin-screw extruder kneading element with two discs, are conducted. They show the optimal convergence behavior of the presented method. In Section 5, we apply our method to compute the temperature-dependent flow inside rotating screw machines. These test cases demonstrate the high potential of the method for complex moving time-dependent simulations of realistic industrial applications. Finally, we draw conclusions in Section 6.

% !TEX root = ./main.tex
\section{A Sliding Mesh Approach for Non-Isothermal Fluid Flow} \label{sec:slidingMesh}

The main motivation for the present work are situations where the domain of the simulation consists of individual parts that share a common interface. Furthermore, the individual parts can rotate or move such that the domain boundaries slide over each other at the common interface.
Very often, the individual domains show complex geometric features.
In case of flow simulation, it is extremely important to accurately represent these geometric features, e.g., small gaps, with a proper discretization.
Examples are single or twin-screw extruders (SSE/TSE) that consist of individual screw parts stacked behind each other.
The resulting individual domains are characterized by extremely small gaps between screw and outer barrel.
In order to capture the flow effects -- especially large gradients in pressure -- correctly, it is of utter importance to resolve the small gaps.
Thus, boundary-conforming meshes can be beneficial.
Special boundary-conforming mesh update techniques have been developed for twin-screw extruders in \cite{helmig2019boundary,hinz2019boundary}.

\begin{figure}[h]
  \centering
  \def\svgwidth{210pt}
  \subfigure[3D view.]{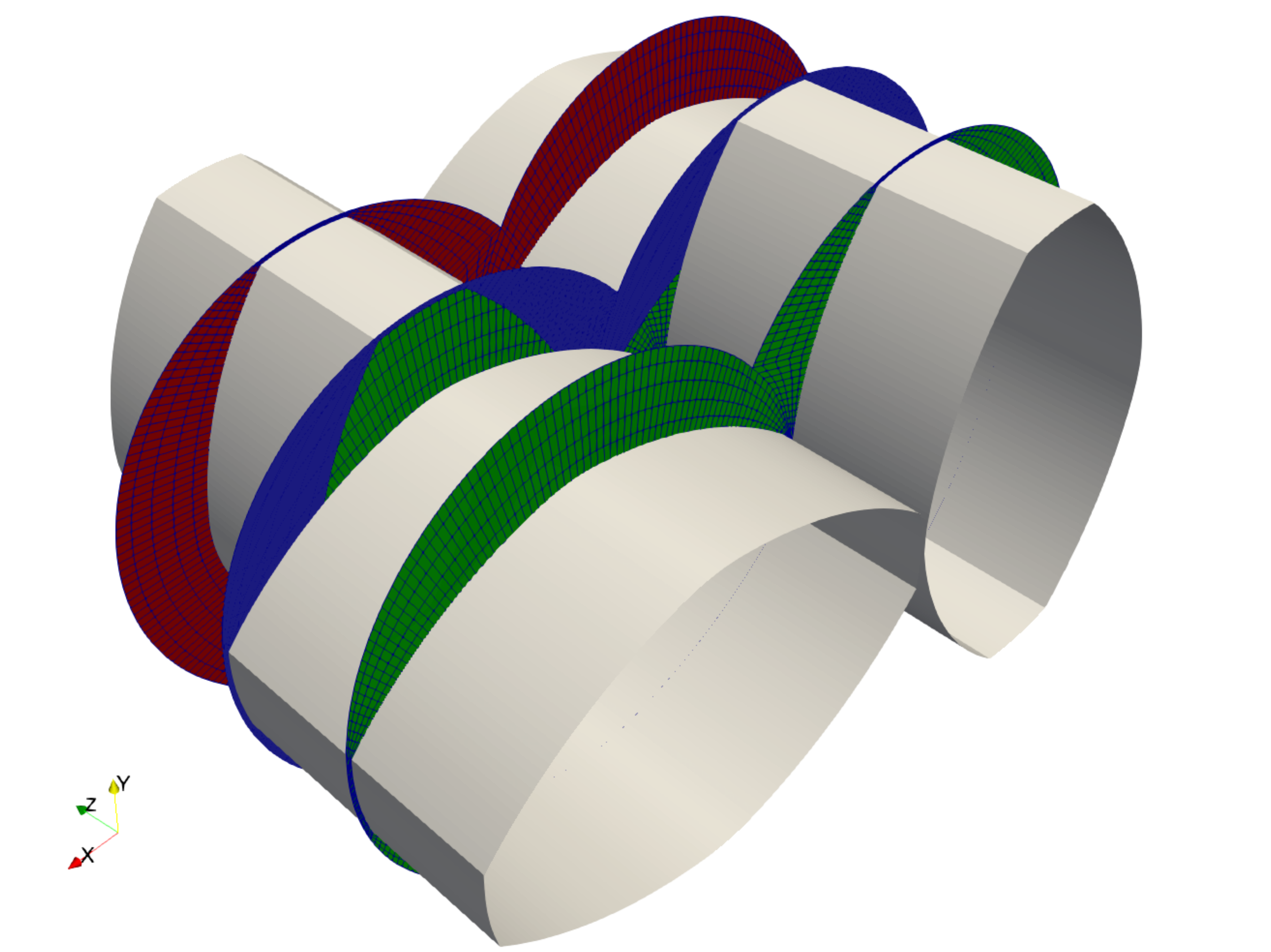\label{fig:slidingMeshExtruder3D}}
  \centering
  \def\svgwidth{240pt}
  \subfigure[2D plane at common interface.]{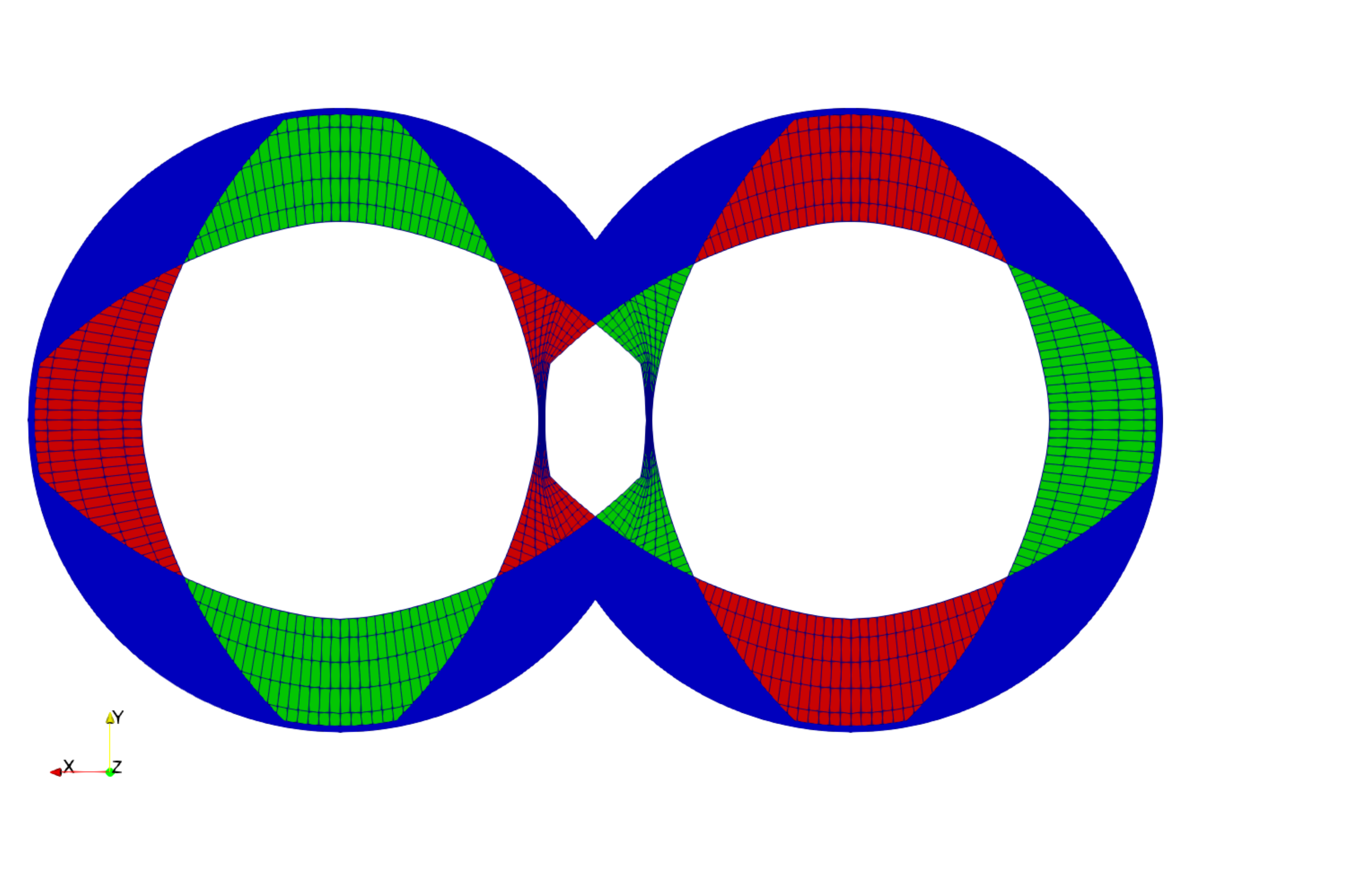\label{fig:slidingMeshExtruder2D}}
  \caption{3D example of twin-screw extruder kneading element with two discs: The screw positions experience a jump of 90$^\circ$ in longitudinal direction across the two discs. The screw of each disc is discretized with a structured boundary-conforming mesh. They are coupled at the common interface $\Gamma_{ff}^{1,2}$ (blue) where the jump occurs.}
  \label{fig:slidingMeshExtruder}
\end{figure}

Fig. \ref{fig:slidingMeshExtruder} shows the example of a so-called kneading element as it can be found in a twin-screw extruder. This kneading element consists of two fixed discs that are staggered by $90^{\circ}$, see Fig. \ref{fig:slidingMeshExtruder3D}.
Structured boundary-conforming meshes can be designed for the individual discs, shown in green and red. However, they don't match at the common interface $\Gamma_{ff}^{1,2}$. Fig. \ref{fig:slidingMeshExtruder2D} shows the two individual meshes at the common interface as well as the overlapping part in blue.
Within the next sections, we will in detail describe a method that enables us to couple the individual domains at non-matching interfaces.

\subsection{Domain Coupling Using a Sliding Interface}

\begin{figure}[h]
    \centering
    \def\svgwidth{250pt}
    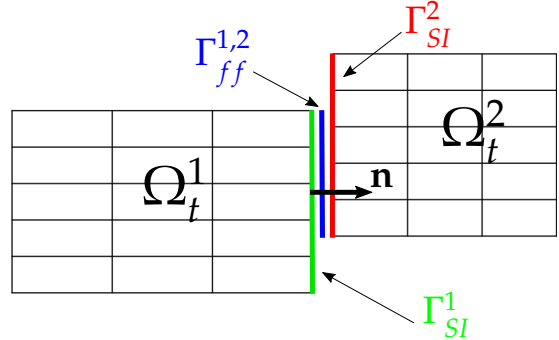
    \caption{Sketch of the sliding mesh approach in 2D.}
    \label{fig:slidingMesh2D}
\end{figure}

We consider a setup in which the overall time-dependent fluid domain $\Omega_t$ is divided into individual spatial subdomains $\Omega^i_t, \; i=1,...,n_{d}$, where $n_d$ is the number of subdomains, see Fig. \ref{fig:slidingMesh2D}.
It holds that $\Omega_t = \bigcup_{i=1}^{n_d} \Omega_t^i$ and $ \bigcap_{i=1}^{n_d} \Omega_t^i = \emptyset$.
Individual subomains $\Omega_t^i$ and $\Omega_t^j$ share a common interface $\Gamma_{ff}^{i,j}$ that can change over time due to a movement of the individual subdomains.
The normal vector to $\Gamma_{ff}^{i,j}$ is defined as ${\bf n} \coloneqq {\bf n}^i = -{\bf n}^j$, where ${\bf n}^i$ and ${\bf n}^j$ are the outward pointing normal vectors of the corresponding boundaries of the subdomains $\Omega^i_t$ and $\Omega^j_t$, respectively.
The sliding boundaries of each subdomain related to $\Gamma_{ff}^{i,j}$ are denoted by $\Gamma_{SI}^i$ and $ \Gamma_{SI}^j$.
Each subdomain $\Omega^i_t$ is discretized using a finite element mesh $\mathcal{T}^i$. We use a classical continous finite element approximation space for each subdomain

\begin{align} \label{eq:functionspaces}
\mathcal{V}^{i,h} := \{ v \in H^1 \left( \Omega^i_t \right) : v |_{K^i} \in \mathcal{Q}^p \left( K^i \right) \forall K^i \in \mathcal{T} ^i \}.
\end{align}

$K^i$ is a tetrahedral or hexahedral element of $\mathcal{T}^i$. $\mathcal{Q}^p$ denotes the polynomials of order $p$ in each direction on $K^i$. In principle the polynomial order can vary on each subdomain. However, only linear polynomials $\mathcal{Q}^1$ will be used within this work.

\subsection{Governing Equations and Coupling Conditions}

We model the fluid as a viscous, incompressible fluid on the moving domain $\Omega _t \; \subset \; \mathbb{R}^{n_{sd}}$, with $n_{sd}$ being the spatial dimension.
As already mentioned, $\Omega _t$ consists of $n_{d}$ disjoint subdomains $\Omega^i_t$ which are enclosed by their boundaries $\Gamma^i _t$, where $t \in (0,T)$ is an instance of time.
The velocity ${\bf u}^i$, pressure $p^i$ and temperature $T^i$ in every subdomain $\Omega^i_t$ are governed by the incompressible Navier-Stokes and heat equations in convective form:

\begin{align}
  \nabla \cdot {\bf{u}}^i = 0 \quad \mbox{on} \ \Omega^i_t, \quad \forall t \in (0,T), \label{eq:cont} \\
 \rho \dfrac{ \hat{\partial} {\bf{u}}^i }{ \partial t } + \rho \left( {\bf{u}}^i - {\bf{u}}_{ALE}^i \right) \cdot \nabla {\bf{u}}^i
 - \nabla \cdot \boldsymbol{\sigma} \left( {\bf{u}}^i, p^i \right) - \rho {\bf b}  = \boldsymbol{0} \quad
 \mbox{on} \ \Omega^i_t, \quad \forall t \in (0,T), \label{eq:momentum} \\
 \rho c_p  \dfrac{ \hat{\partial} T^i}{ \partial t} + \rho c_p \left( \boldsymbol{u}^i - \boldsymbol{u}_{ALE}^i \right) \cdot \nabla T^i  - \kappa \boldsymbol {\Delta} T^i - 2 \eta \nabla \boldsymbol{ u}^i \colon \boldsymbol{\varepsilon} \left( \boldsymbol{u}^i \right)  \; = 0 \quad
  \mbox{on} \ \Omega^i_t, \quad \forall t \in (0,T), \label{eq:heat}
\end{align}

where $\rho$ is the fluid density, ${\bf b}$ the gravity vector, $\kappa$ the thermal conductivity and $c_p$ the specific heat capacity.
$\dfrac{ \hat{\partial} (\cdot) }{ \partial t }$ represents the time derivative in the arbitrary Lagrangian Eulerian (ALE) frame and ${\bf{u}}_{ALE}^i$ is the domain mesh velocity. A detailed derivation of the ALE description can be found in \cite{donea2003finite,forster2006geometric}. The Cauchy stress tensor $\boldsymbol{\sigma}$ is defined as:

\begin{align}
    \boldsymbol{\sigma} ( {\boldsymbol{u} }^i, p^i) = -p^i {\boldsymbol{I} } + 2 \eta \left( \dot{\gamma }, T^i  \right) \boldsymbol{\varepsilon}( {\boldsymbol{u} }^i), \\
    \boldsymbol{\varepsilon}( {\boldsymbol{u} }^i) = \frac{1}{2} \left( \nabla \boldsymbol{ u}^i + \left( \nabla \boldsymbol{ u}^i  \right)^T \right),
\end{align}

 with $\boldsymbol{\varepsilon}( {\boldsymbol{u} }^i)$ being the strain-rate tensor and $\eta$ the dynamic viscosity. $\eta$ is constant for Newtonian fluids, and for Generalized Newtonian models varies with respect to temperature $T$ and shear rate $\dot{\gamma}$. \\

The latter is defined as
 %The Generalized Newtonian models are used in case we aim at modeling the flow of plastic melt to account for shear-thinning behavior. For these models the viscosity depends on the invariants of the rate of strain tensor $\boldsymbol{\varepsilon}$, such as the shear rate $\dot{\gamma}$:

 \begin{align}
     \dot{\gamma} = \sqrt{2 \boldsymbol{\varepsilon} \left( {\bf u}^i \right)  \colon \boldsymbol{\varepsilon} \left( {\bf u}^i \right) }.
 \end{align}

 Within this work we use two different models, namely the Carreau and the Cross model with WLF correction.\\

 The Carreau model \cite{carreau1979review} states:

 \begin{align}
     \eta \left( \dot{\gamma} \right) = \eta _{\infty}  + \left( \eta_0 - \eta _{\infty} \right)  \left( 1+ \left( \lambda \dot{\gamma} \right) ^2 \right)^{\frac{n-1}{2}},
 \end{align}

 where $\lambda$ is the relaxation time, $n$ is the power index, $\eta_0 $ is the viscosity at zero shear rate and $\eta _{\infty} $ is the viscosity at infinite shear rate. \\

 The Cross model \cite{rudolph2014polymer} is:

 \begin{align}
     \eta \left( \dot{\gamma} \right) = \frac{\eta_0 }{1+\left( \frac{\eta_0  \dot{\gamma}}{\tau ^*} \right)^{(1-n)}},
 \end{align}
 %\begin{align}
%     \eta \left( \dot{\gamma}, T^i \right) = \frac{\eta_0 \left( T^i \right)}{1+\left( \frac{\eta_0  \left( T^i \right) \dot{\gamma}}{\tau ^*} \right)^{(1-n)}},
 %\end{align}

 where $\tau ^*$ is the critical shear stress at the transition from the Newtonian plateau.

Furthermore, we also want to model the influence of temperature on the viscosity. This is done by making the viscosity at zero shear rate $\eta_0 $ dependent on temperature via the WLF correction:

 \begin{align}
     \eta_{0} (T^i) = D_1 \; exp \left( - \frac{A_1  \left( T^i - T_{ref} \right) }{A_2 + \left( T^i - T_{ref} \right) } \right),
 \end{align}

 where $D_1$ is the viscosity at a reference temperature $T_{ref}$ and $A_1$ and $A_2$ are parameters that describe the temperature dependency. \\

The Dirichlet and Neumann boundary conditions for flow and temperature fields are:

\begin{align}
    {\bf u}^i = {\bf g}^{i,f} \; \mbox{on} \; \left( \Gamma_t^i \right) ^{f} _g, \\
    \boldsymbol{\sigma} ( {\boldsymbol{u} }^i, p^i) \cdot {\bf n} = {\bf h}^{i,f} \; \mbox{on} \; \left( \Gamma_t^i \right) ^{f} _h, \\
    T^i = g^{i,t} \; \mbox{on} \; \left( \Gamma_t^i \right) ^{t} _g, \\
    \kappa \nabla T^i \cdot {\bf n} = h^{i,t} \; \mbox{on} \; \left( \Gamma_t^i \right) ^{t} _h,
\end{align}

where superscripts $f$ and $t$ denote flow and temperature, respectively, and
$\left( \Gamma_t^i \right) ^{f,t} _g$ and $\left( \Gamma_t^i \right) ^{f,t} _h$ are portions of $\left(\Gamma _t^i \right) ^{f,t}$. \\

Furthermore, we need to conserve mass, momentum and energy over the common interface $\left(\Gamma_t^{i,j}\right)_{ff}$ between two subdomains $\Omega_t^i$ and $\Omega_t^j$ . The mass conservation in combination with a noslip condition requires continuity of velocity and temperature:

\begin{align}\label{eq:jumpmass}
    \llbracket {\bf u} \rrbracket = {\bf 0} \quad \mbox{on} \quad \left( \Gamma _t^{i,j} \right)_{ff}, \nonumber \\
    \llbracket T \rrbracket = 0 \quad \mbox{on} \quad \left( \Gamma_t^{i,j} \right)_{ff}.
\end{align}

Momentum and energy conservation are obtain by demanding equal surface tractions and heat fluxes:

\begin{align}\label{eq:jumpflux}
    \llbracket \boldsymbol{\sigma} ( {\boldsymbol{u} }, p)] \cdot {\bf n} \rrbracket = {\bf 0} \quad \mbox{on} \quad \left( \Gamma _t^{i,j} \right)_{ff}, \nonumber \\
    \llbracket \kappa \nabla T \cdot {\bf n} \rrbracket = 0 \quad \mbox{on} \quad \left( \Gamma _t^{i,j} \right)_{ff}.
\end{align}

$\llbracket \cdot \rrbracket$ denotes the jump in a function over the interface $ \left( \Gamma_t \right)^{i,j}_{ff}$ and is defined as $\llbracket \cdot \rrbracket = (\cdot)^i - (\cdot)^j$.

\section{Stabilized Finite Element Discretization for Sliding Meshes} \label{sec:stabFEMSlidingMesh}

\subsection{Finite Element Discretization}

In the following, we discretize the equations \eqref{eq:cont} -- \eqref{eq:heat}. We use a SUPG/PSPG type stabilized formulation for each subdomain $\Omega_t^i$ following \cite{pauli2017stabilized}. It can also be interpreted as a variational multiscale formulation as presented in \cite{hughes2018multiscale}.
The weak form is derived by multiplying the continuity equation \eqref{eq:cont}, momentum equation \eqref{eq:momentum} and heat equation \eqref{eq:heat} with test functions $ q^{i}$, ${\bf w}^i$ and $ v^{i}$ respectively.
Next, we choose appropriate finite dimensional test and trial functions ${\bf w}^{i,h}$, $q^{i,h}$, $v^{i,h}$, ${\bf u}^{i,h}$, $p^{i,h}$, $T^{i,h}$ from \eqref{eq:functionspaces} and integrate by parts taking the boundary conditions into account.
We define the $L^2$-inner product on $\Omega_t^i$ as $(\cdot,\cdot)_{\Omega_t^i}$ and as $( \cdot,\cdot )_{\left( \Gamma _t^{i,j} \right)_{ff}}$  on $\left( \Gamma _t^{i,j} \right)_{ff}$.
In the following, we use the simplified notation $(\cdot,\cdot)_{\Omega_t} =  \sum_{i=1}^{n_{d}} (\cdot,\cdot)_{\Omega_t^i}$ and
$( \cdot,\cdot )_{\Gamma_{ff}} = \sum_{i=1}^{n_{d}} \sum_{j=i+1}^{n_{d}}( \cdot,\cdot )_{\left( \Gamma _t^{i,j} \right)_{ff}}$.
Adding the stabilization terms, we end up with the stabilized weak formulation for the fluid equations

\begin{align}\label{eq:weakflow}
    \mbox{B}^f \left( \left[ {\bf w}^h, q^h \right]; \left[ {\bf u}^h, p^h \right] \right) + \mbox{I}^f \left( \left[ {\bf w}^h, q^h \right]; \left[ {\bf u}^h, p^h \right] \right) = 0,
\end{align}

where

\begin{align}\label{eq:stabgalerkinflow}
    \mbox{B}^f \left( \left[ {\bf w}^h, q^h \right]; \left[ {\bf u}^h, p^h \right] \right) = \;
    &\left( {\bf w}^h , \rho \frac{ \hat{\partial} {\bf u}^h }{\partial t} \right)_{\Omega_t} + \left( {\bf w}^h, \rho ( {\bf u}^h - {\bf u}_{ALE}^h ) \cdot \nabla {\bf u }^h  \right)_{\Omega_t} \\
    &+ \left( q^h, \nabla \cdot {\bf u}^h \right)_{\Omega_t} - \left( \nabla \cdot {\bf w}^h, p^h \right)_{\Omega_t} + \left( \boldsymbol{\varepsilon} ({\bf w}^h), 2 \eta \boldsymbol{\varepsilon} ({\bf u}^h) \right)_{\Omega_t} \nonumber \\
    &+ \left( \rho ( \frac{ \hat{\partial} {\bf w}^h }{\partial t} + ( {\bf u}^h - {\bf u}_{ALE}^h ) \cdot \nabla {\bf w }^h), \tau_{\mbox{\tiny{MOM}}} {\bf r}_{\mbox{\tiny{MOM}}}^h \right)_{\Omega_t} \nonumber \\
    &+ \left( \nabla \cdot {\bf w}^h, \tau_{\mbox{\tiny{CON}}} r_{\mbox{\tiny{CON}}}^h \right)_{\Omega_t}
    + \left( \nabla q^h, \tau_{\mbox{\tiny{MOM}}} {\bf r}_{\mbox{\tiny{MOM}}}^h \right)_{\Omega_t} \nonumber \\
    &-  \left( {\bf w}^h, \rho {\bf b} \right)_{\Omega_t}
    - \left( {\bf w}^h, {\bf h}^f \right) _{\Gamma^f_h}, \nonumber
\end{align}

and

\begin{align}\label{eq:jumpfluxflow}
    \mbox{I}^f \left( \left[ {\bf w}^h, q^h \right]; \left[ {\bf u}^h, p^h \right] \right) =
    \left( \llbracket {\bf w}^h, p^h {\bf n} - 2 \eta \boldsymbol{\varepsilon} ({\bf u}^h) ) \cdot {\bf n} \rrbracket \right) _{\Gamma_{ff}}.
\end{align}

The stabilized weak formulation for the temperature equation is

\begin{align} \label{eq:weaktemp}
    \mbox{B}^t \left( v^h ; T^h  \right) + \mbox{I}^t \left( v^h ; T^h  \right) = 0,
\end{align}

where

\begin{align}\label{eq:stabgalerkintemp}
    \mbox{B}^t \left( v^h ; T^h  \right) = & \left( v^h , \rho c_p \frac{ \hat{\partial} T^h }{\partial t} \right)_{\Omega_t} + \left( v^h, \rho c_p ( {\bf u}^h - {\bf u}_{ALE}^h ) \cdot \nabla T^h  \right)_{\Omega_t} \\
    & + \left( \nabla v^h, \kappa \nabla T^h \right)_{\Omega_t} \nonumber \\
    %&- \left( \llbracket v^h, \kappa \nabla T^h \cdot {\bf n} \rrbracket \right) _{\Gamma_{ff}} \\
    &+ \left( \rho c_p ( \frac{ \hat{\partial} v^h }{\partial t} + ( {\bf u}^h - {\bf u}_{ALE}^h ) \cdot \nabla v^h), \tau_{\mbox{\tiny{TEMP}}} r_{\mbox{\tiny{TEMP}}}^h \right)_{\Omega_t} \nonumber \\
    & - \left( v^h, 2 \eta \nabla {\bf u}^h \colon \boldsymbol{\varepsilon} ( {\bf u}^h ) \right)_{\Omega_t}
    - \left( v^h, h^t \right) _{\Gamma^t_h} \nonumber
\end{align}

and
\begin{align}\label{eq:jumpfluxtemp}
    \mbox{I}^t \left( v^h ; T^h  \right) = - \left( \llbracket v^h, \kappa \nabla T^h \cdot {\bf n} \rrbracket \right) _{\Gamma_{ff}}.
\end{align}

$\mbox{I}^f \left( \left[ {\bf w}^h, q^h \right]; \left[ {\bf u}^h, p^h \right] \right)$ and $\mbox{I}^t \left( v^h ;  T^h \right)$
are the interface consistency or jump flux terms.
Note that we did not include the enforcement of the interface coupling condition yet. This will be done in the next section. \\
The first line of equations \eqref{eq:stabgalerkinflow} and \eqref{eq:stabgalerkintemp} are the standard Galerkin terms.
The remaining lines of equations \eqref{eq:stabgalerkinflow} and \eqref{eq:stabgalerkintemp} are the residual-based stabilization terms,
where the residuals $r_{\mbox{\tiny{CON}}}^h$, ${\bf r}_{\mbox{\tiny{MOM}}}^h$ and $r_{\mbox{\tiny{TEMP}}}^h$ are defined as:

\begin{align}
    &r_{\mbox{\tiny{CON}}}^h = \nabla \cdot {\bf{u}}^h, \\
    &{\bf r}_{\mbox{\tiny{MOM}}}^h = \rho \dfrac{ \hat{\partial} {\bf{u}}^h }{ \partial t }
    + \rho \left( {\bf{u}}^h - {\bf{u}}_{ALE}^h \right) \cdot \nabla {\bf{u}}^h
    - \nabla \cdot \boldsymbol{\sigma} \left( {\bf{u}}^h, p^h \right) - \rho {\bf b}^h, \label{eq:residualmomentum}\\
    &r_{\mbox{\tiny{TEMP}}}^h = \rho c_p  \dfrac{ \hat{\partial} T^h}{ \partial t} +
    \rho c_p \left( \boldsymbol{u}^h - \boldsymbol{u}_{ALE}^h \right) \cdot \nabla T^h
    - \kappa \boldsymbol {\Delta} T^h - 2 \eta \nabla \boldsymbol{ u}^h \colon \boldsymbol{\varepsilon} \left( \boldsymbol{u}^h \right). \label{eq:residualheat}
\end{align}

The stress contributions in equations \eqref{eq:residualmomentum} and \eqref{eq:residualheat} involve second order derivatives. In case only first order polynomials are used, we recover these terms by using a least-squares recovery technique \citep{jansen1999better}.
This improves the consistency of the stabilized method especially for highly viscous problems.
The terms in line three of equation \eqref{eq:stabgalerkinflow} as well as of equation \eqref{eq:stabgalerkintemp} are the SUPG stabilization terms.
They are used to stabilize the formulations for convection-dominated problems.
The first term in line four of equation \eqref{eq:stabgalerkinflow} adds artificial diffusion to stabilize the continuity equation.
The second term is the PSPG term which is needed to make the formulation inf-sub stable, since we use equal order polynomial for all degrees of freedom.
The stabilization parameters $\tau_{\mbox{\tiny{MOM}}}$, $\tau_{{\mbox{\tiny{CONT}}}}$ and $\tau_{\mbox{\tiny{TEMP}}}$ are based on expressions given in \citep{pauli2017stabilized,pauli2016stabilized}.

\subsection{Nitsche Coupling at Common Interface} \label{sec:nitscheCoupling}

Following \cite{schott2016stabilized}, we can now formulate the Nitsche coupling. First, we have to define the weighted interface average operator:

\begin{align}\label{eq:weightedinterface}
    \{a \} \coloneqq k^i a^i + k^j a^j, \quad \langle a \rangle \coloneqq k^j a^i + k^i a^j \quad \mbox{with} \quad k^i, k^j \geq 0 \quad \mbox{and} \quad k^i + k^j = 1.
\end{align}

$k^i$ and $k^j$ are real positive weights. Typical choices are one-sided weightings, meaning $k^i=1, k^j=0$, or an averaged weighting also denoted as $\{ \cdot \}_m$ where $k^i=k^j=0.5$ as used, e.g., in \cite{bazilevs2008nurbs}. However, this might lead to an unbalanced weighting between strains in case generalized Newtonian models are used. This is due to possible large jumps in viscosity between elements. Thus, we will use a weighting that balances the difference in viscosity, given as: $k^i = \eta^j /(\eta^i + \eta^j), \;k^j = \eta^i /(\eta^i + \eta^j)$. We use the same weighting for the heat equation but substitute $\eta$ by $\kappa$.
In case of constant value for $\eta $ or  $ \kappa$, it results in the standard averaged weighting.

Using the relations given in \eqref{eq:weightedinterface}, we can reformulate the interface jump terms \eqref{eq:jumpfluxflow} and \eqref{eq:jumpfluxtemp} under the assumption that

\begin{align}
    \llbracket ab \rrbracket = \langle a \rangle \llbracket b \rrbracket + \llbracket a \rrbracket \{ b \}
\end{align}

and by incorporating the flux condition \eqref{eq:jumpflux} to:

\begin{align}
    \left( \llbracket {\bf w}^h, p^h \cdot {\bf n} - 2 \eta \boldsymbol{\varepsilon} ({\bf u}^h) ) \cdot {\bf n} \rrbracket \right) _{\Gamma_{ff}} &=
    \left( \llbracket  {\bf w}^h \rrbracket , \{ p^h \} \cdot {\bf n} \right) _{\Gamma_{ff}}
    - \left( \llbracket {\bf w}^h \rrbracket , \{ 2 \eta \boldsymbol{\varepsilon} ({\bf u}^h)  \cdot {\bf n} \} \right) _{\Gamma_{ff}}, \label{eq:interfacejumpflow}\\
    - \; \left( \llbracket v^h, \kappa \nabla T^h \cdot {\bf n} \rrbracket \right) _{\Gamma_{ff}} &=
    - \; \left( \llbracket  v^h \rrbracket , \{ \kappa \nabla T^h \cdot {\bf n} \}  \right) _{\Gamma_{ff}}. \label{eq:interfacejumptemp}
\end{align}

Using a Nitsche coupling inspired by \cite{schott2016stabilized,bazilevs2008nurbs} we can extend the interface jump term for the fluid $\mbox{I}^f \left( \left[ {\bf w}^h, q^h \right]; \left[ {\bf u}^h, p^h \right] \right)$ (see Eq. \eqref{eq:interfacejumpflow}), and obtain the following stabilized formulation for the flow coupling condition at the common interface $\Gamma_{ff}$:

\begin{align} \label{eq:nitscheflow}
    \mbox{I}^f_{stab} \left( \left[ {\bf w}^h, q^h \right]; \left[ {\bf u}^h, p^h \right] \right) = & \; \left( \llbracket  {\bf w}^h \rrbracket , \{ p^h \} \cdot {\bf n} \right) _{\Gamma_{ff}}
    - \left( \llbracket {\bf w}^h \rrbracket , \{ 2 \eta \boldsymbol{\varepsilon} ({\bf u}^h)  \cdot {\bf n} \} \right) _{\Gamma_{ff}} \nonumber \\
    &- \; \left(  \{ q^h \} \cdot {\bf n}  , \llbracket {\bf u}^h \rrbracket \right) _{\Gamma_{ff}}
    - \left(  \{ 2 \eta \boldsymbol{\varepsilon} ({\bf w}^h)  \cdot {\bf n} \}, \llbracket {\bf u}^h \rrbracket \right) _{\Gamma_{ff}} \nonumber \\
    &+ \; \left( \tau_{\mbox{\tiny{SI}}}^{f} \llbracket {\bf w}^h \rrbracket, \llbracket {\bf u}^h \rrbracket \right) _{\Gamma_{ff}} \\
    &+ \; \left( \rho ( ({\bf u}^h - {\bf u}^h_{ALE}) \cdot {\bf n} ) \{ {\bf w}^h \}_m, \llbracket {\bf u}^h \rrbracket \right) _{\Gamma_{ff}}
    - \frac{1}{2} \left( \rho | ({\bf u}^h - {\bf u}^h_{ALE}) \cdot {\bf n} | \llbracket {\bf w}^h \rrbracket, \llbracket {\bf u}^h \rrbracket \right) _{\Gamma_{ff}}. \nonumber
\end{align}

The terms in the second line are the so-called adjoint-consistency terms.
They fulfill the coupling condition \eqref{eq:jumpmass} and ensure mass conservation over the interface.
Note that the terms involving the pressure trial (first line) and test function (second line) are skew-symmetric.
This ensures that the form is stability-neutral without violating the adjoint-consistency of the formulation. The second to last line is a consistent penalization term that ensures coercivity of the formulation and furthermore, ensures mass conservation over the interface.
We define $\tau_{\mbox{\tiny{SI}}}^{f}$ as:

\begin{align} \label{eq:stabnitscheflow}
    \tau_{\mbox{\tiny{SI}}}^{f} = \frac{\alpha}{2} \frac{\eta^i \eta^j}{\eta^i + \eta^j} \left( \frac{1}{h^i} + \frac{1}{h^j} \right),
\end{align}

with

\begin{align}
    h^i = 2 \left(  ({\bf n}^i)^T \cdot \hat{{\bf G}}^i \cdot {\bf n}^i \right)^{-\frac{1}{2}},
    \quad \hat{{\bf G}} = \left( \frac{\partial \boldsymbol{\xi} }{\partial {\bf x} } \right)^T \left( \frac{\partial \boldsymbol{\xi} }{\partial {\bf x} } \right).
\end{align}

$\frac{\partial \boldsymbol{\xi}}{\partial {\bf x} }$ is the inverse Jacobian of the element mapping between reference and physical domain and $\alpha$ is a stabilization parameter. The last line in equation \eqref{eq:nitscheflow} is an upwinding stabilization that controls instabilities when mass is transported from one subdomain into another one \cite{schott2016stabilized}. \\

The coupling for the temperature at $\Gamma_{ff}$ follows the same Nitsche approach already applied to the flow. Thus, similarly to Eq. \eqref{eq:nitscheflow}, we extend $\mbox{I}^t \left( v^h; T^h \right) $ ( see Eq. \eqref{eq:interfacejumptemp}) and obtain:

\begin{align} \label{eq:nitschetemp}
    \mbox{I}^t_{stab} \left( v^h; T^h \right) = &- \; \left( \llbracket  v^h \rrbracket , \{ \kappa \nabla T^h \cdot {\bf n} \}  \right) _{\Gamma_{ff}}
    - \left( \{ \kappa \nabla v^h \cdot {\bf n} \} , \llbracket  T^h \rrbracket  \right) _{\Gamma_{ff}}
    + \left( \tau_{\mbox{\tiny{SI}}}^{t} \llbracket v^h \rrbracket, \llbracket T^h \rrbracket \right) _{\Gamma_{ff}} \nonumber \\
    &+ \; \left( \rho c_p ( ({\bf u}^h - {\bf u}^h_{ALE}) \cdot {\bf n} ) \{ v^h \}_m, \llbracket T^h \rrbracket \right) _{\Gamma_{ff}}
    - \frac{1}{2} \left( \rho c_p | ({\bf u}^h - {\bf u}^h_{ALE}) \cdot {\bf n} | \llbracket v^h \rrbracket, \llbracket T^h \rrbracket \right) _{\Gamma_{ff}}.
\end{align}

We define $\tau_{\mbox{\tiny{SI}}}^{t}$ similar to \eqref{eq:stabnitscheflow} as:

\begin{align} \label{eq:stabnitschetemp}
    \tau_{\mbox{\tiny{SI}}}^{t} = \frac{\alpha}{2} \frac{\kappa^i \kappa^j}{\kappa^i + \kappa^j} \left( \frac{1}{h^i} + \frac{1}{h^j} \right).
\end{align}

\subsection{Weak Imposition of Boundary Conditions Using a Multimesh Technique}

Within this section, we will discuss the case in which the interface between two domains does not fully match geometrically, i.e., $\Gamma_{SI} \neq \Gamma_{ff}^i \cap \Gamma_{ff}^j$. In this case, we only want the fluid to flow through the overlap, see Fig. \ref{fig:slidingMesh2D}.
In the remaining part of the boundary interface we apply a no-slip condition to impose the movement of the underlying rigid parts of the other subdomain.
Thus, we need to be able to apply a Dirichlet boundary condition in the part of $\Gamma_{SI}$ that is not connected to the other subdomain, namely $\Gamma_{SI} \backslash \Gamma^{i,j}_{ff}$.
The Dirichlet boundary conditions will be imposed weakly also using Nitsche's method \cite{nitsche1971variationsprinzip,bazilevs2007weak}. This is done by adding the following terms to the stabilized weak form of the flow equations:

\begin{align} \label{eq:weakdcflow}
    \mbox{B}^{f} &\left( \left[ {\bf w}^h, q^h \right]; \left[ {\bf u}^h, p^h \right] \right) + \mbox{I}^{f}_{stab} \left( \left[ {\bf w}^h, q^h \right]; \left[ {\bf u}^h, p^h \right] \right) \nonumber \\
    &+ \; \left( {\bf w}^h, p^h {\bf n} \right) _{\Gamma_{SI} \backslash \Gamma^{i,j}_{ff}}
    - \left( {\bf w}^h, 2 \eta \boldsymbol{\varepsilon} ({\bf u}^h)  \cdot {\bf n} \right) _{\Gamma_{SI} \backslash \Gamma^{i,j}_{ff}} \nonumber \\
    &- \; \left( q^h {\bf n}, ( {\bf u}^h - {\bf u}^h_g ) \right) _{\Gamma_{SI} \backslash \Gamma^{i,j}_{ff}}
    - \left( 2 \eta \boldsymbol{\varepsilon} ({\bf w}^h)  \cdot {\bf n}, ( {\bf u}^h - {\bf u}^h_g ) \right) _{\Gamma_{SI} \backslash \Gamma^{i,j}_{ff}} \\
    &+ \left( {\bf w}^h , \alpha \eta / h \; ( {\bf u}^h - {\bf u}^h_g ) \right) _{\Gamma_{SI} \backslash \Gamma^{i,j}_{ff}} = 0 \nonumber.
\end{align}

The same concept is also applied to weakly impose Dirichlet temperature boundary conditions on $\Gamma_{SI} \backslash \Gamma^{i,j}_{ff}$:

\begin{align} \label{eq:weakdctemp}
    \mbox{B}^{t} &\left(  v^h ;  T^h \right) + \mbox{I}^{t}_{stab} \left( v^h ; T^h \right) \nonumber \\
    &- \; \left( v^h , \kappa \nabla T^h \cdot {\bf n} \right) _{\Gamma_{SI} \backslash \Gamma^{i,j}_{ff}}
    - \left( \kappa \nabla v^h \cdot {\bf n}, (T^h - T^h_g)  \right) _{\Gamma_{SI} \backslash \Gamma^{i,j}_{ff}} \\
    &+ \; \left( v^h, \alpha \kappa / h \; (T^h  - T^h_g) \right) _{\Gamma_{SI} \backslash \Gamma^{i,j}_{ff}} = 0. \nonumber
\end{align}

\subsection{Algorithmic Details}

In order to evaluate integrals on $\Gamma^{i,j}_{ff}$ as they appear in \eqref{eq:nitscheflow} and \eqref{eq:nitschetemp} it is necessary to compute the cuts between boundary elements $K^i$ on $\Gamma_{SI}^i$ and $K^j$ on $\Gamma_{SI}^j$ in order to define quadrature rules. In contrast to embedded mesh techniques where the dimension of the cut cells is equal to the dimension of the domain, the computation of cut cells for our approach is one dimension lower than that of the computational domain. Thus, we only need to compute cut cells in 2D. The computation of the cut cells is performed using CGAL \cite{cgal:eb-20a}.
First, we detect all possible collisions of boundary elements of the subdomains at all common sliding interfaces using bounding boxes.
The candidate cut elements are then checked for collision using geometric predicates.
In the following, we compute the intersection of two elements $K^i$ and $K^j$.
$K^i \cap K^j$ is convex such that we can simply triangulate the cut polygon to define standard Gaussian quadrature rules on the resulting sub-triangles.

\begin{figure}[h]
    \centering
    \def\svgwidth{250pt}
    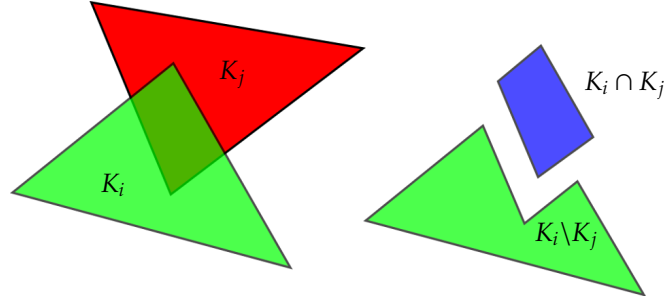
    \caption{Sketch of the multimesh approach: two intersecting triangles (green and red) and the resulting cut polygon (blue).}
    \label{fig:multiMesh}
\end{figure}

The open question is how to determine the quadrature rule for the weak imposition of Dirichlet boundary conditions on $\Gamma_{SI} \backslash \Gamma^{i,j}_{ff}$.
Defining quadrature rules for all elements on the boundary that are not cut is trivial.
However, we also have to take the partly cut elements into account.
The uncut part of an element is not necessarily convex, which makes sub-triangulation non-trivial (see Fig. \ref{fig:multiMesh}).
Thus, we will make use of the multimesh concept of nested quadrature introduced in \cite{johansson2019multimesh}.
Considering Fig. \ref{fig:multiMesh}, we can interpret $K^i$ as an element on $\Gamma_{SI}^i$ and $K^j$ as an element on $\Gamma_{SI}^j$ such that $K^i \cap K^j$ is the resulting cut element. We aim to integrate on the uncut part of $K^i$ denoted as $K^i \backslash K^j$.
As described before, we are able to define quadrature rules on $K^i \cap K^j$ since we can present it as a set of triangles.
The idea of Johansson {\it et al} \cite{johansson2019multimesh} is to not explicitly define a quadrature rule for $K^i \backslash K^j$. Instead they use the inclusion-exclusion principle of combinatorics:

\begin{align}
    | K^i \cup K^j | = |K^i| + |K^j| - |K^i \cap K^j|
\end{align}

to derive

\begin{align}
    |K^i \backslash K^j | &= |K^i \cup K^j| \backslash |K^j| \nonumber \\
    & = |K^i| + |K^j| - |K^i \cap K^j| - |K^j| \\
    & = |K^i| - |K^i \cap K^j|. \nonumber
\end{align}

Thus, we can simply integrate on $K^i \backslash K^j$ by integrating on $K^i$ and also on $K^i \cap K^j$ but using negative weights.
Extending this concept to the integration of weak Dirichlet conditions on the sliding interface, we can simply integrate the terms on $\Gamma_{SI} \backslash \Gamma^{i,j}_{ff}$ by integrating on $\Gamma_{SI}$ using standard quadrature rules and then simply integrate the terms on all cut elements on $\Gamma_{ff}$ using negative weights.

\begin{remark}
 Summing up and neglecting terms everywhere on the boundary due to the negative quadrature weights sounds like a computationally expensive process at first. However, the terms of the weak imposition of boundary condition \eqref{eq:weakflow} and those for the Nitsche coupling at the common interface \eqref{eq:nitscheflow} have a very similar structure. Thus, basically no extra computational cost due to the subtracting arises during assembly, since the terms are already computed for the Nitsche coupling.
\end{remark}

% !TEX root = ./main.tex
%!TeX spellcheck = en-US
\section{Numerical Examples} \label{sec:numericalexamples}

\subsection{2D Taylor-Green Flow}

We study the convergence behavior of our method computing the 2D Taylor-Green problem \cite{pearson1964computational} for a convection-dominated ($\eta=0.0001$, $Re=10000$) and viscous-dominated ($\eta=0.1$, $Re=10$) flow regime. The analytical solution for the flow in a domain $\Omega = [0,1]^2$ is given as:

\begin{align} \label{eq:analsoltaylorgreen}
    &{\bf u} = \left( -sin(2\pi y)cos(2\pi x) e^{-8 \pi^2 \eta t} ,  sin(2\pi x)cos(2\pi y) e^{-8 \pi^2 \eta t} \right), \nonumber \\
    &p       = - \frac{1}{4} \left( cos(4\pi x) + cos(4\pi y)  \right)  e^{-16 \pi^2 \eta t}.
\end{align}

\begin{figure}[h]
  \centering
  \subfigure[Pressure field for refinement level 3.]{\includegraphics[width=.48\linewidth]{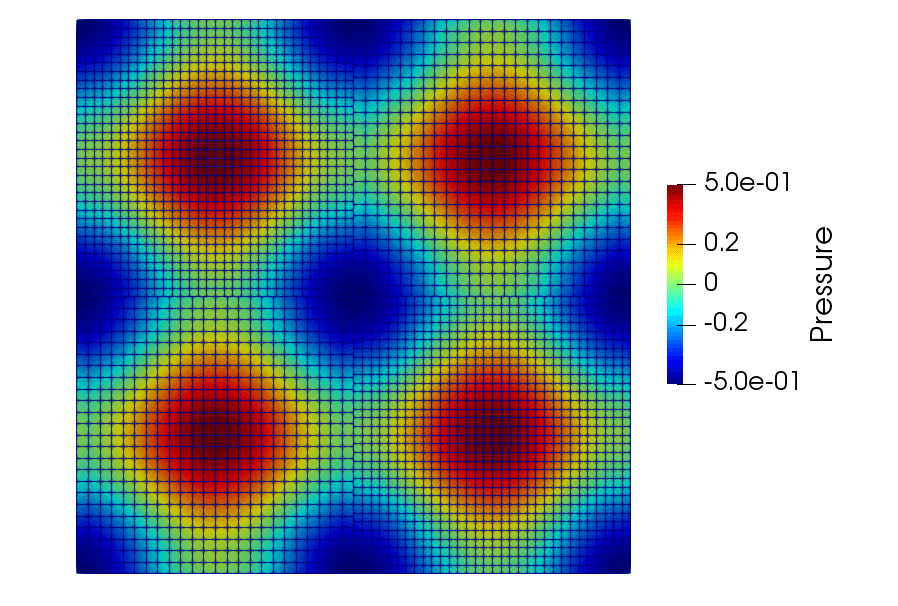}\label{fig:pressure2DTaylorGreen}}
  \centering
  \subfigure[Velocity magnitude contour on mesh with refinement level 2.]{\includegraphics[width=.48\linewidth]{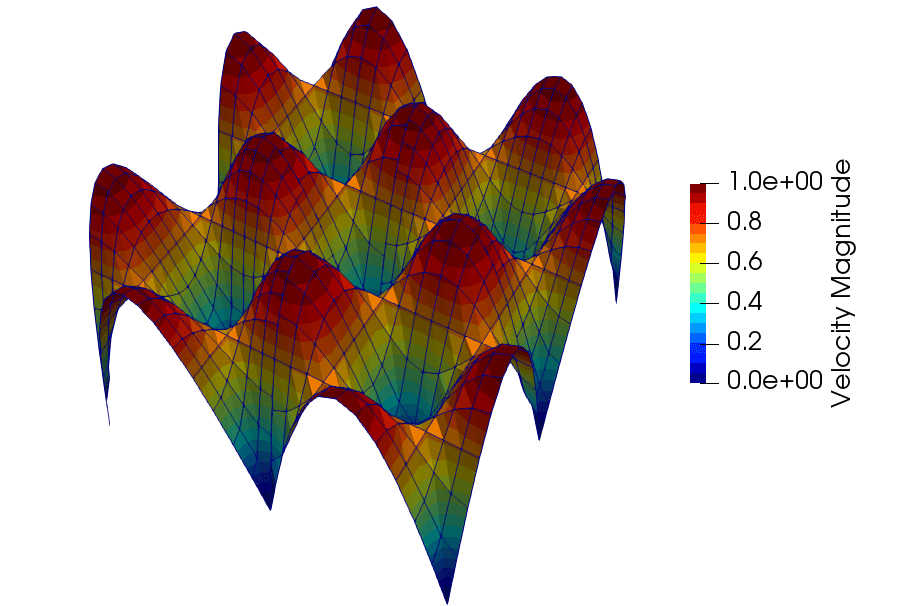}\label{fig:velocity2DTaylorGreen}}
  \caption{Velocity and pressure fields for stationary Taylor-Green flow with $\eta = 0.1$.}
  \label{fig:result2DTaylorGreen}
\end{figure}

This problem has already been used to test approaches based on Nitsche's method for embedded as well as cut mesh methods \cite{schott2016stabilized,dokken2019multimesh,massing2018stabilized}.
We  divide the domain into four equally sized square domains (see Fig. \ref{fig:pressure2DTaylorGreen}). The two subdomains on the upper left and lower right are discretized with a mesh with $n$ elements in each direction, and the off-diagonal ones with a mesh with $m$ elements in each direction.
For the convergence study, we start with $n=4$ and $m=3$ and refine seven times.
In order to compare our results, we also use one mesh that is not subdivided.
The initial number of elements for the reference mesh in one direction is 8.
We set the analytical solution as Dirichlet boundary condition on the outer boundary.
In order to have a constant pressure level, we set the analytical solution for the pressure in the lower left corner.
For the viscous case with $\eta = 0.1$, we compute a steady solution.
It is obtained by setting the body force components ${\bf b}$ to the negative time derivative of the analytical flow solution.
The steady solution for pressure and velocity is shown in Fig. \ref{fig:result2DTaylorGreen}.
For the convection-dominated case, we compute 10 time steps with time step size $\Delta t = 0.00025$. We discretize the time derivative using a Backwards Differentiation Formula of first order (BDF1) \cite{forti2015semi}.
The very small time step is used in order to prevent that time discretization errors affect the convergence study.
As initial condition we set the analytical flow solution (Eq. \eqref{eq:analsoltaylorgreen}) for $t = 0$.
We analyze the L$^2$ error norm on the whole domain

\begin{align}
    \left\lVert {\bf u}^h - {\bf u} \right\rVert _{L^2(\Omega)} \quad \mbox{and} \quad \left\lVert p^h - p \right\rVert _{L^2(\Omega)},
\end{align}

as well as the interface coupling error norm

\begin{align} \label{eq:interfaceErrorNorm}
    \left\lVert \llbracket {\bf u}^h \rrbracket \right\rVert _{L^2(\Gamma^{ff})} \quad \mbox{and} \quad \left\lVert \llbracket p^h \rrbracket \right\rVert _{L^2(\Gamma^{ff})}.
\end{align}

The error norms for the velocity should converge with second order and the pressure errors norms with first order \cite{burman2014fictitious,massing2014stabilized}.

\subsubsection{Spatial Convergence Study} \label{sec:spatialConvergence2D}

\begin{figure}
  \centering
  \subfigure[]{\includegraphics[width=.45\linewidth]{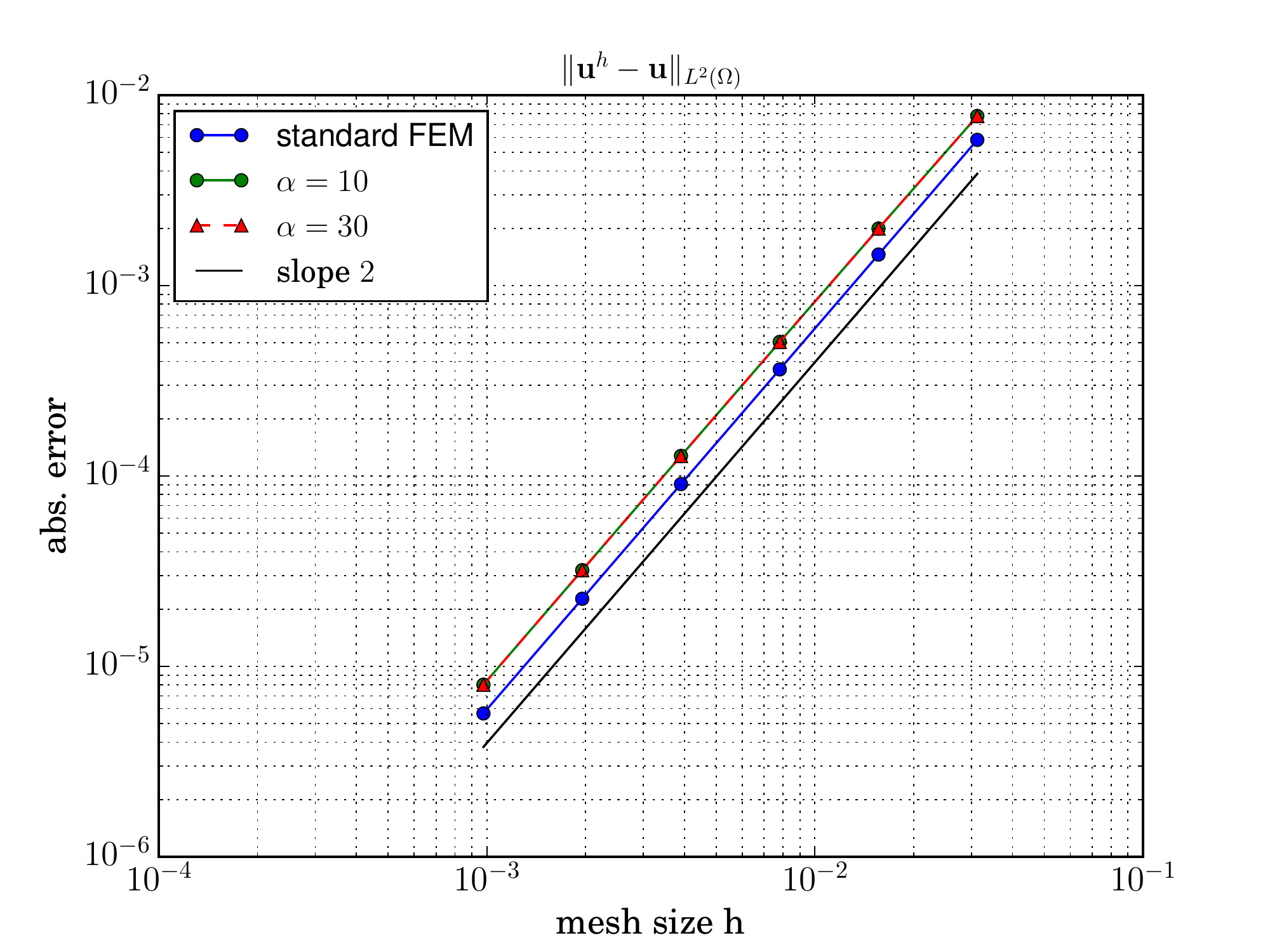}\label{fig:convergence2DTaylorGreenConvectiveVel}}
  \centering
  \subfigure[]{\includegraphics[width=.45\linewidth]{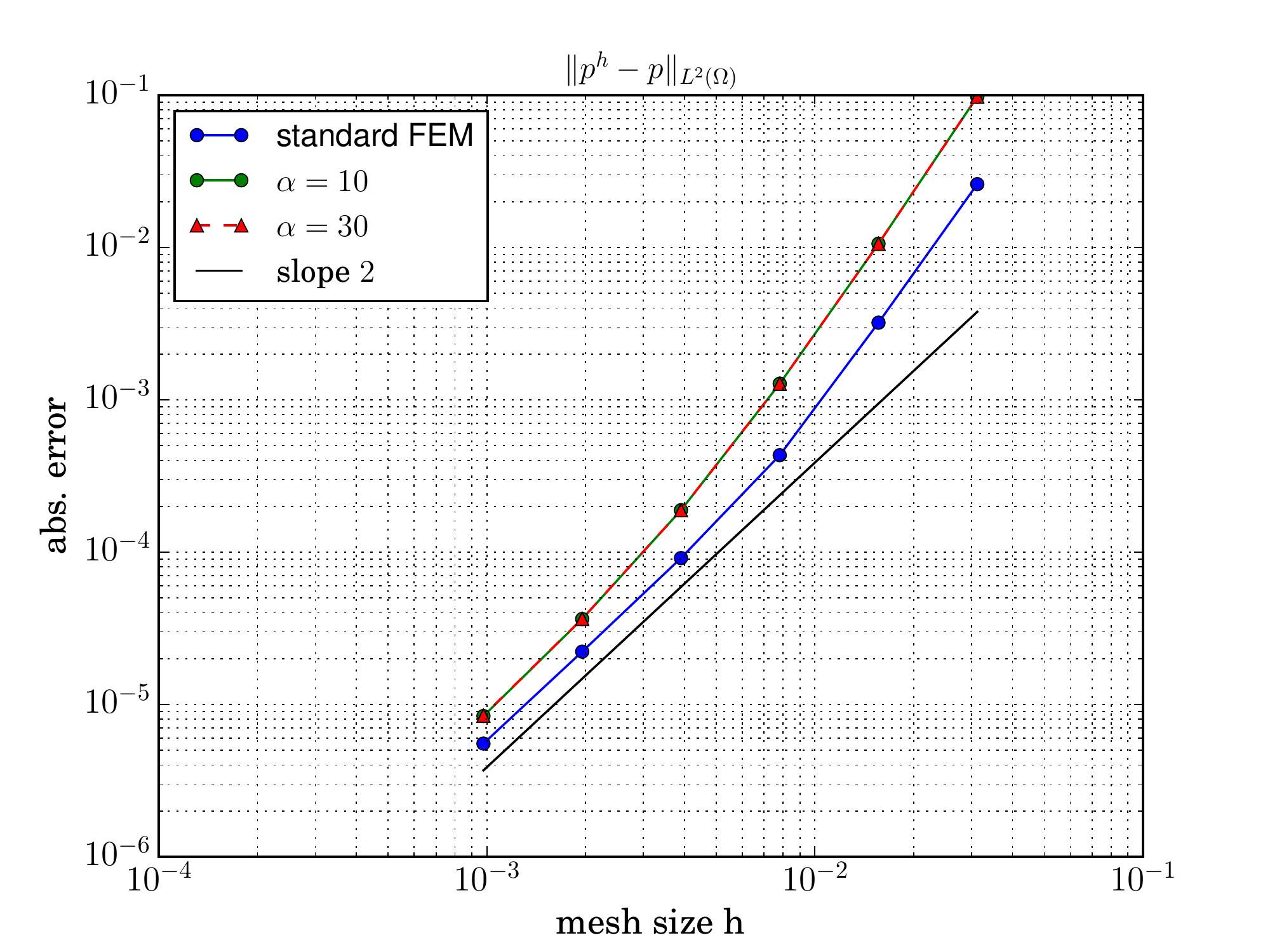}\label{fig:convergence2DTaylorGreenConvectivePressure}}
  \centering
  \subfigure[]{\includegraphics[width=.45\linewidth]{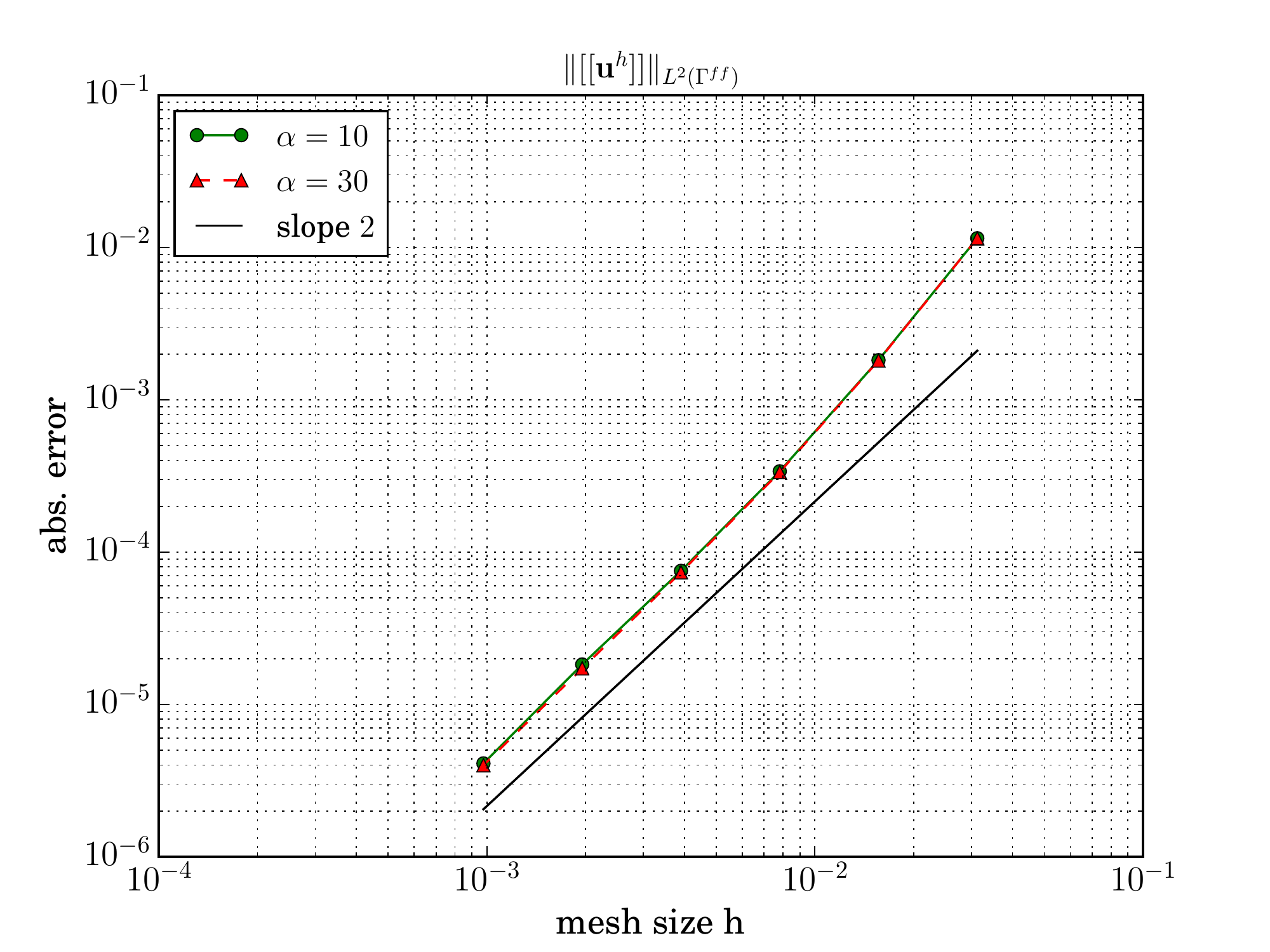}\label{fig:convergence2DTaylorGreenConvectiveVelInterface}}
  \centering
  \subfigure[]{\includegraphics[width=.45\linewidth]{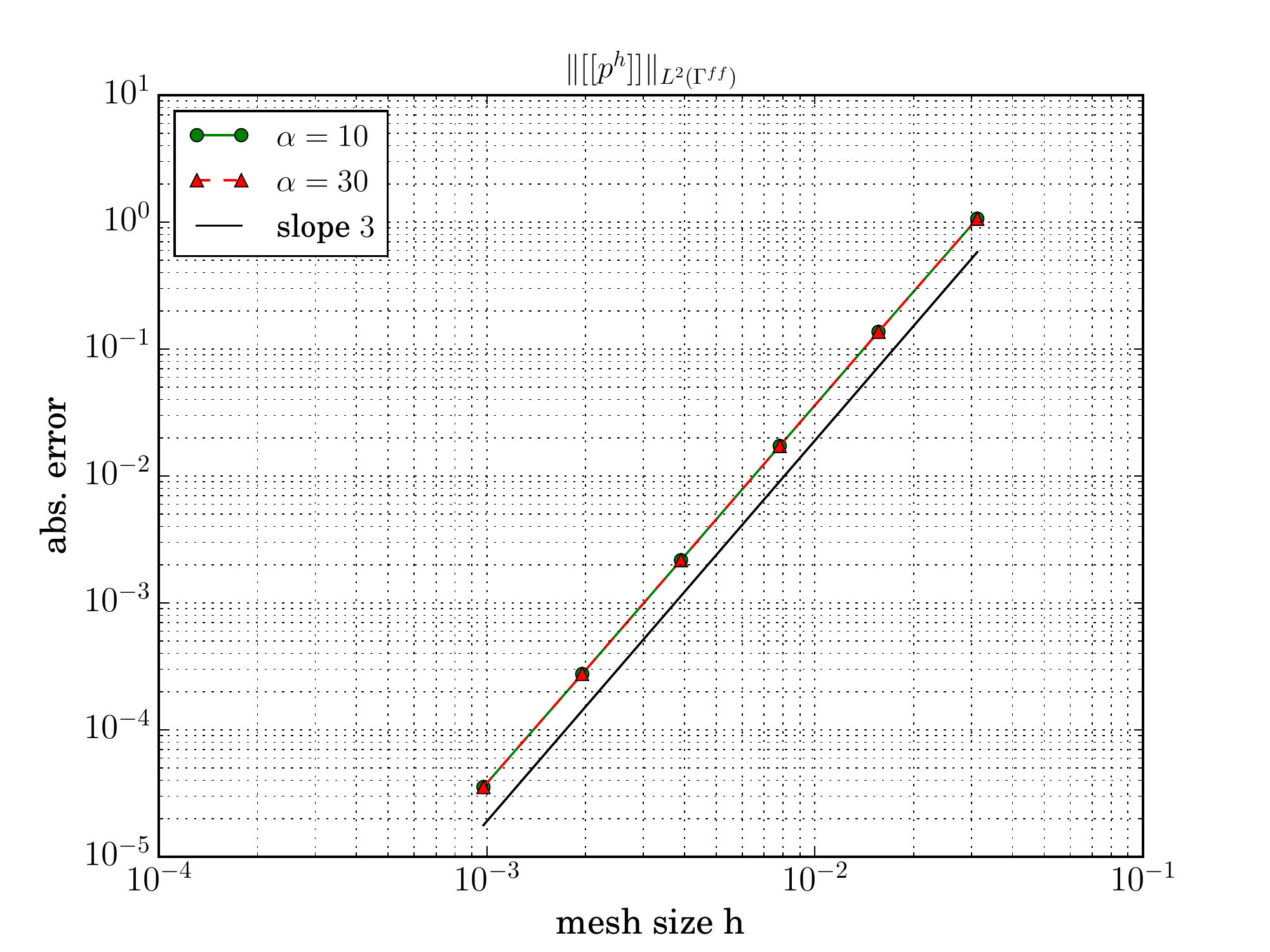}\label{fig:convergence2DTaylorGreenConvectivePressureInterface}}
  \caption{Spatial convergence for Taylor-Green flow with $\eta = 0.0001$: domain error norms (top row) and interface error norms (bottom row). Label 'standard FEM' refers to standard finite elements on an uncut mesh and $\alpha=10,30$ is the sliding mesh approach with the respective value for the Nitsche stabilization factor $\alpha$, see Eq. \eqref{eq:stabnitscheflow}.}
  \label{fig:convergence2DTaylorGreenConvective}
\end{figure}

In a first step, we consider the convergence behavior for convection-dominated flow. The convergence rates are shown in Fig. \ref{fig:convergence2DTaylorGreenConvective}. For the L$^2$ error on the whole domain we compare the results with those on a matching, uncut mesh.
Furthermore, we use two different Nitsche stabilization parameters $\alpha=10$ and $\alpha=30$.
These values are inspired by numerical studies performed in \cite{schott2016stabilized,massing2018stabilized}.

We observe that the domain as well as interface error norms converge with the expected or even higher convergence rates.
The higher rates can be explained by the high regularity of the solution.
Furthermore, Fig. \ref{fig:convergence2DTaylorGreenConvectiveVelInterface} shows that we obtain optimal convergence for the mass conservation across the interface, which indicates that the upwinding scheme as well as the scaling of our stabilization parameters works as expected.
We observe a large difference in the error for the matching and the non-matching case, see Fig. \ref{fig:convergence2DTaylorGreenConvectiveVel} and \ref{fig:convergence2DTaylorGreenConvectivePressure}. This can be explained by the fact that we plot the norm over the minimal element length that is the same for the matching and non-matching case. However, the initial number of elements for the non-matching case in two subdomains is only 3 instead of 4 (refer to Fig. \ref{fig:pressure2DTaylorGreen}). Thus, the overall number of degrees of freedom is always smaller, which explains the offest in the error norm.

\begin{figure}
  \centering
  \subfigure[]{\includegraphics[width=.45\linewidth]{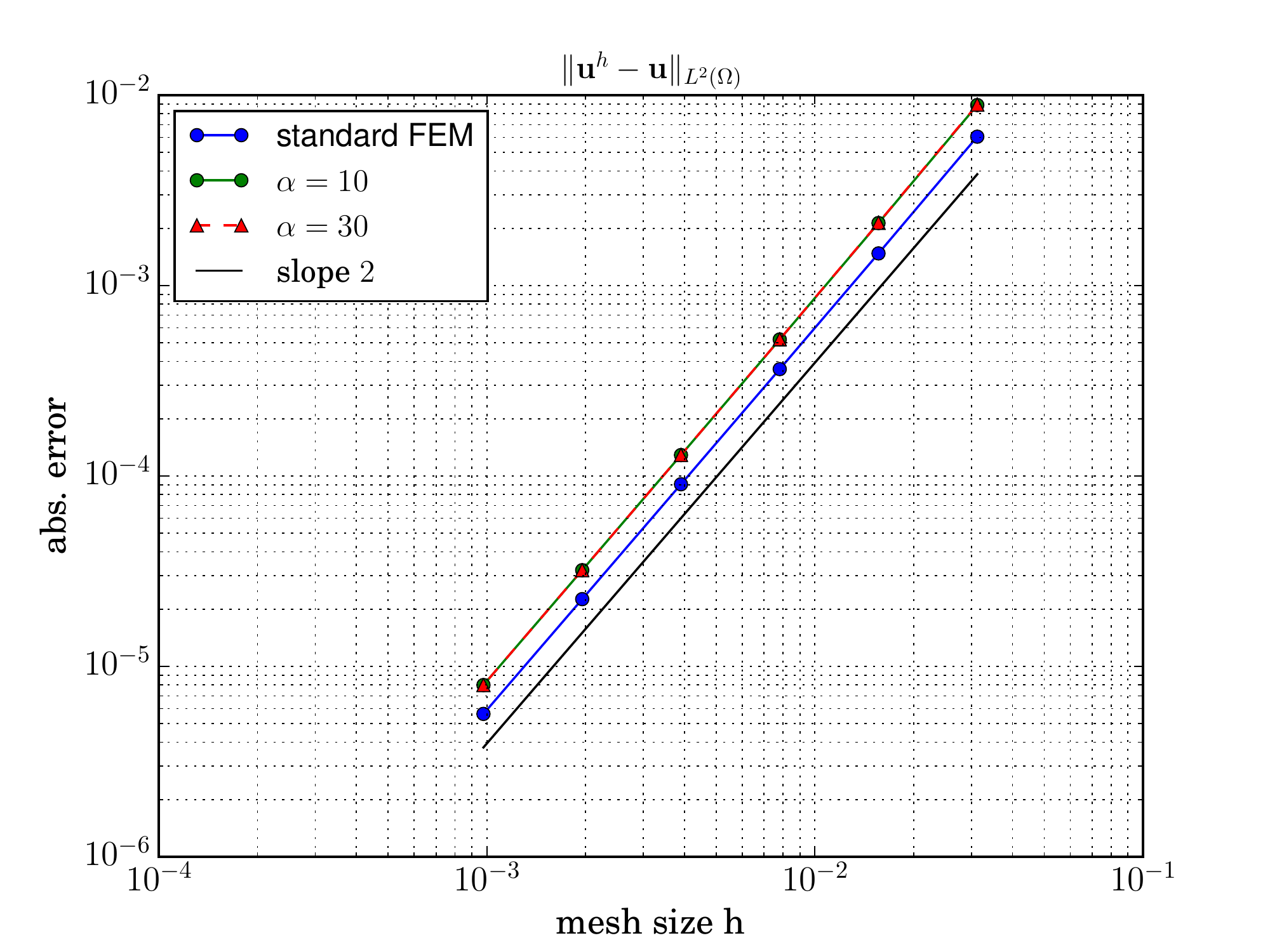}\label{fig:convergence2DTaylorGreenViscousVel}}
  \centering
  \subfigure[]{\includegraphics[width=.45\linewidth]{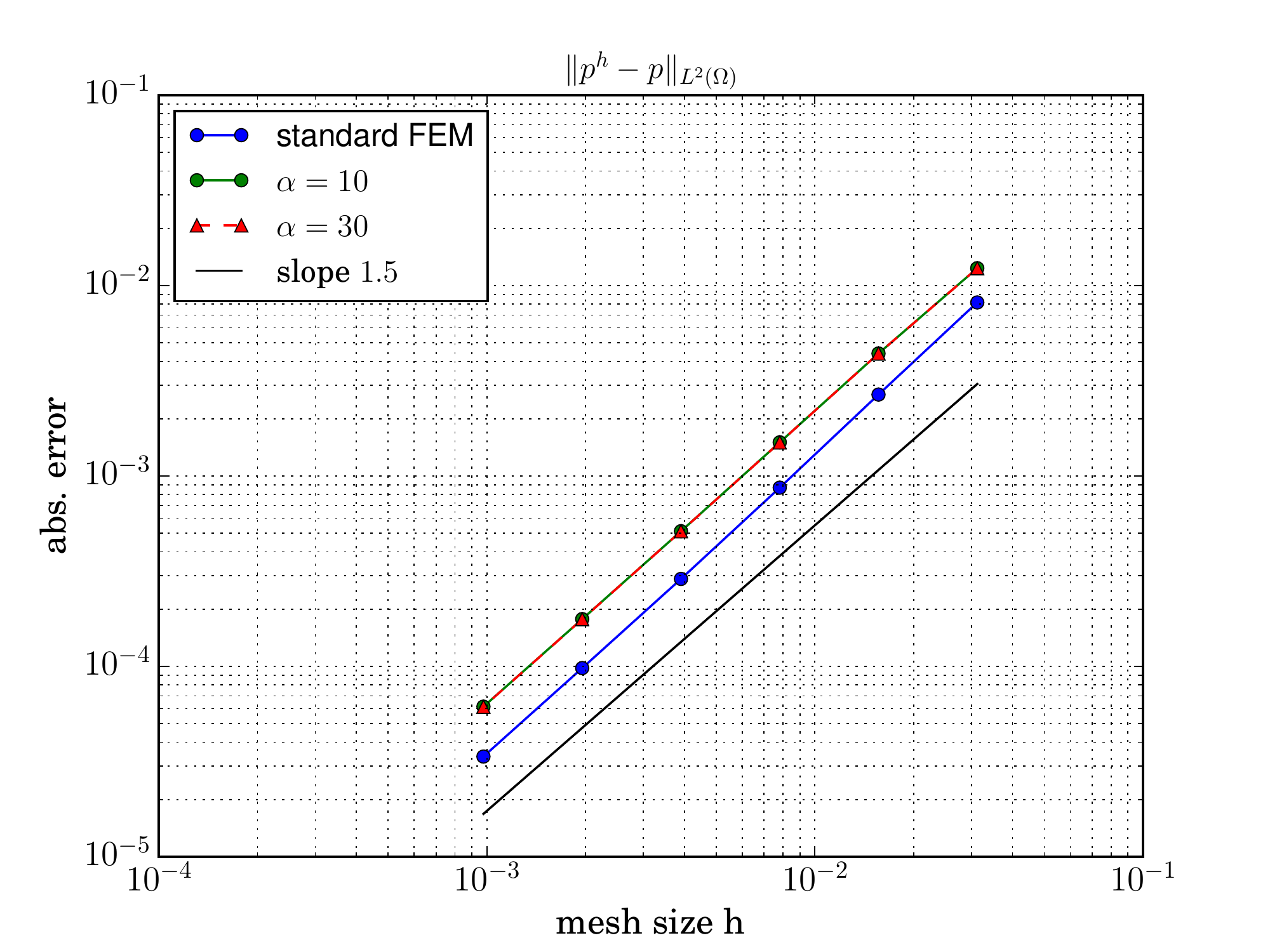}\label{fig:convergence2DTaylorGreenViscousPressure}}
  \centering
  \subfigure[]{\includegraphics[width=.45\linewidth]{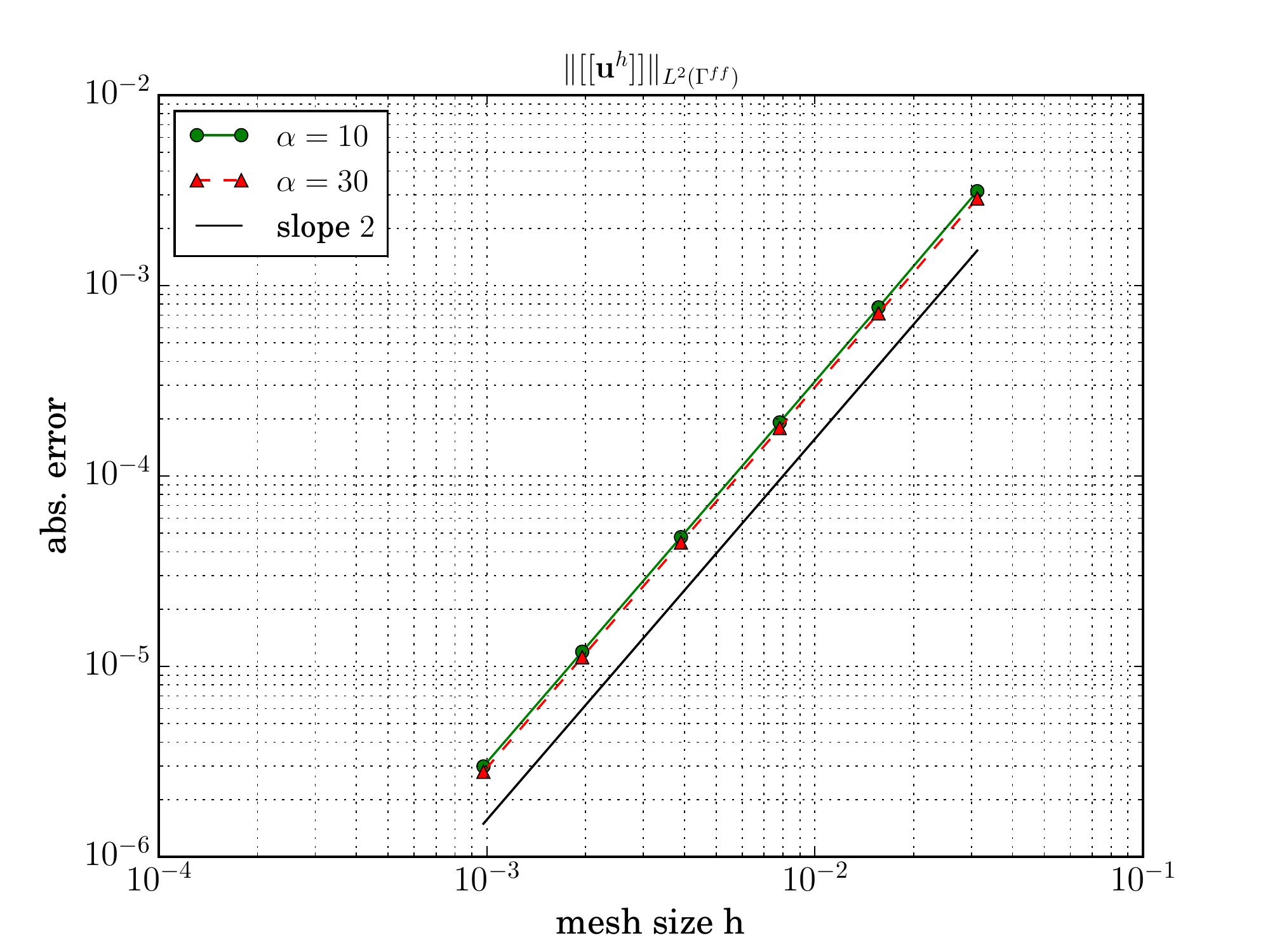}\label{fig:convergence2DTaylorGreenViscousVelInterface}}
  \centering
  \subfigure[]{\includegraphics[width=.45\linewidth]{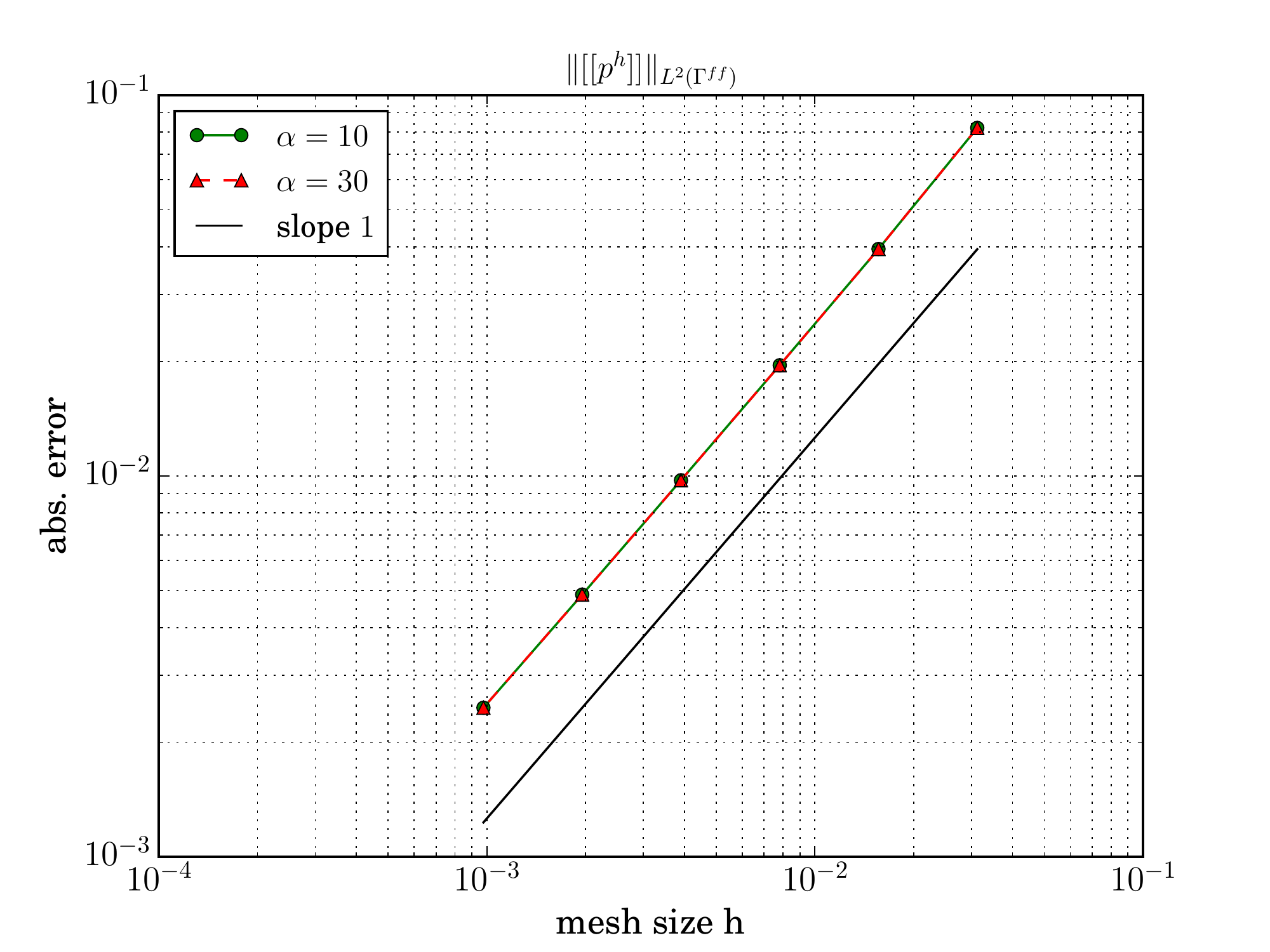}\label{fig:convergence2DTaylorGreenViscousPressureInterface}}
  \caption{Spatial convergence for Taylor-Green flow with $\eta = 0.1 \; Pa \; s$: domain error norms (top row) and interface error norms (bottom row). Label 'standard FEM' refers to standard finite elements on an uncut mesh and $\alpha=10,30$ is the sliding mesh approach with the respective value for the Nitsche stabilization factor $\alpha$, see Eq. \eqref{eq:stabnitscheflow}.}
  \label{fig:convergence2DTaylorGreenViscous}
\end{figure}

We also analyze the spatial convergence for the viscous case.
Fig. \ref{fig:convergence2DTaylorGreenViscous} shows the convergence behavior in all error norms.
Similar to the convective case, we observe optimal convergence rates for the domain as well as interface error norms.
Having a closer look at the error in mass conservation over the interface (see Fig. \ref{fig:convergence2DTaylorGreenViscousVelInterface}), we observe that we obtain slightly smaller errors for $\alpha = 30$ compared to $\alpha = 10$. Thus, $\alpha = 30$ will be used for the following test cases.\\
These optimal spatial convergence rates indicate the correct scaling of our stabilization parameters for different flow regimes. We have not tested the coupling for the temperature. However, the structure of the heat equation is the same as for the momentum equation, so that similar convergence results can be expected.

\subsection{3D Flow in the Kneading Element of a Twin-Screw Extruder}

In this section, we extend the spatial convergence study to a more complex flow problem in 3D.
Specifically, we analyze the flow inside a twin-screw extruder kneading element with two discs of length $20 \;mm$ that are staggered by $90^{\circ}$.
The design of the underlying screw is based on Booy's law \cite{sarhangi2012adaptive,booy1978geometry} and the parameters are given in Table \ref{table:screw3DConvergence}.

\begin{figure}[h]
  \centering
  \includegraphics[width=.6\linewidth]{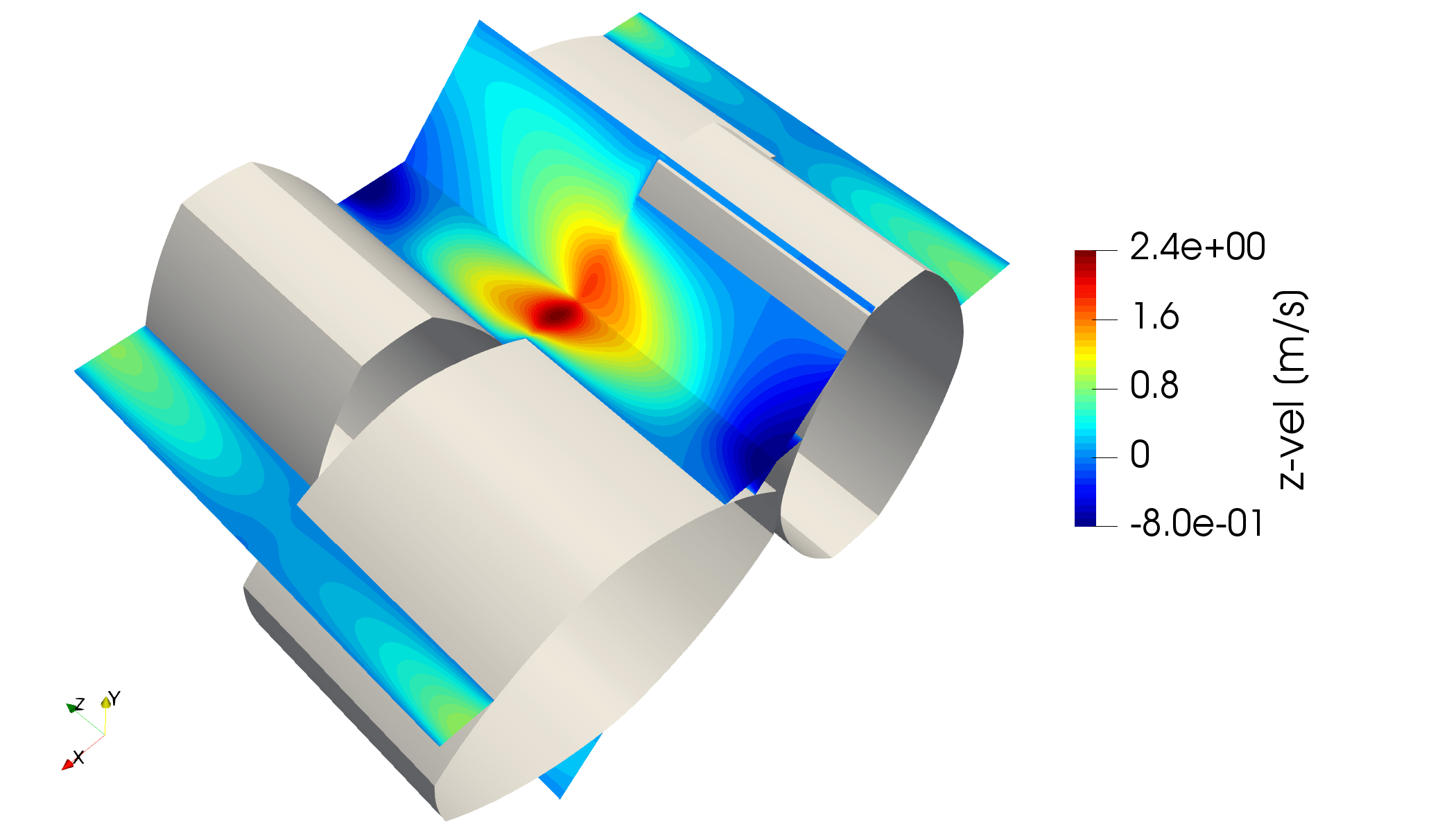}
  \caption{Kneading element with two disc staggered by $90^{\circ}$: the $z$-velocity of the resulting flow field in two planes.}
  \label{fig:extruder3DzVel}
\end{figure}

\begin{table}[h]
\begin{minipage}[t]{0.4\linewidth}
    \centering
    \begin{tabular}{l r}
     \hline
     screw radius $R_s$  & 14.7 $mm$ \\
     center line distance $C_l$ & 26.2 $mm$ \\
     screw-screw clearance $\delta _s$ & 0.3 $mm$ \\
     screw-barrel clearance $\delta _b$ & 0.3 $mm$ \\
     %Pitch length & 20 $mm$ \\
     \hline
     \end{tabular}
     \caption{Geometry parameters for 3D screw element.}
     \label{table:screw3DConvergence}
\end{minipage}
\centering
\begin{minipage}[t]{0.2\linewidth}
    \centering
    \begin{tabular}{l c c}
    \hline
    $\eta_0$  & 1290 & $Pa \; s$ \\
    $\eta_{\infty}$ & 0 & $Pa \; s$ \\
    $n$ & 0.559 & - \\
    $\lambda$ & 0.112 & $s$ \\
    \hline
    \end{tabular}
    \caption{Carreau parameters}
    \label{table:carreau2D}
\end{minipage}
\centering
\begin{minipage}[t]{0.39\linewidth}
    \centering
    \begin{tabular}{l c c c}
     \hline
     mesh & $n_s$  & $n_r$ & $n_a$ \\
     1 & 180 & 5 & 90 \\
     2 & 360 & 10 & 180 \\
     3 & 540 & 15 & 270 \\
     4 & 720 & 20 & 360 \\
     5 & 1080 & 30 & 540 \\
     \hline
     \end{tabular}
     \caption{Mesh discretization based on SRMUM for 3D convergence study. A more detailed description of the parameters can be found in \cite{helmig2019boundary}.}
     \label{table:mesh3DConvergence}
\end{minipage}
\end{table}

\begin{figure}[h]
  \centering
  \subfigure[front view]{\includegraphics[width=.49\linewidth]{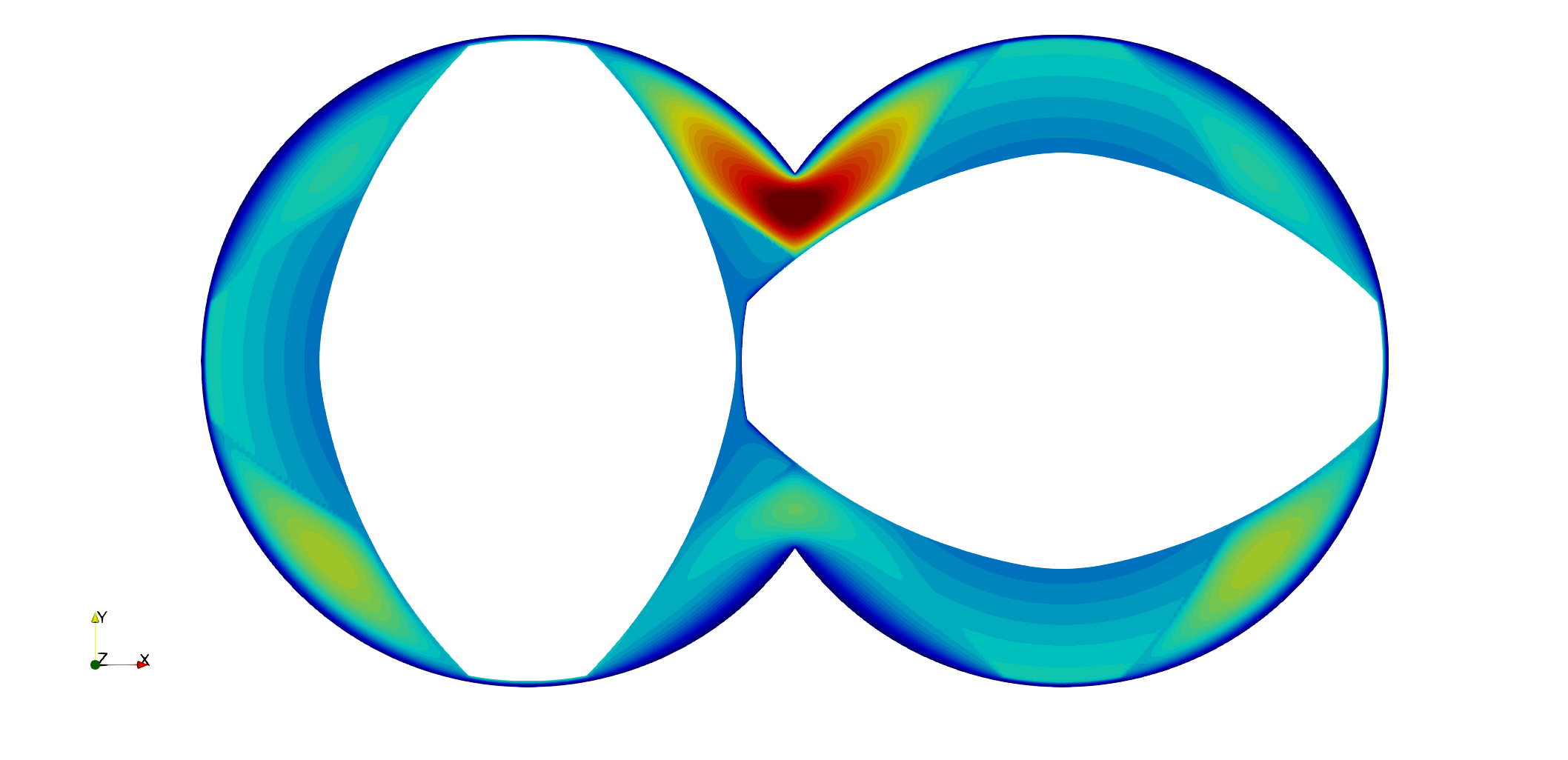}\label{fig:velMag3D2DExtruder1}}
  \centering
  \subfigure[back view]{\includegraphics[width=.49\linewidth]{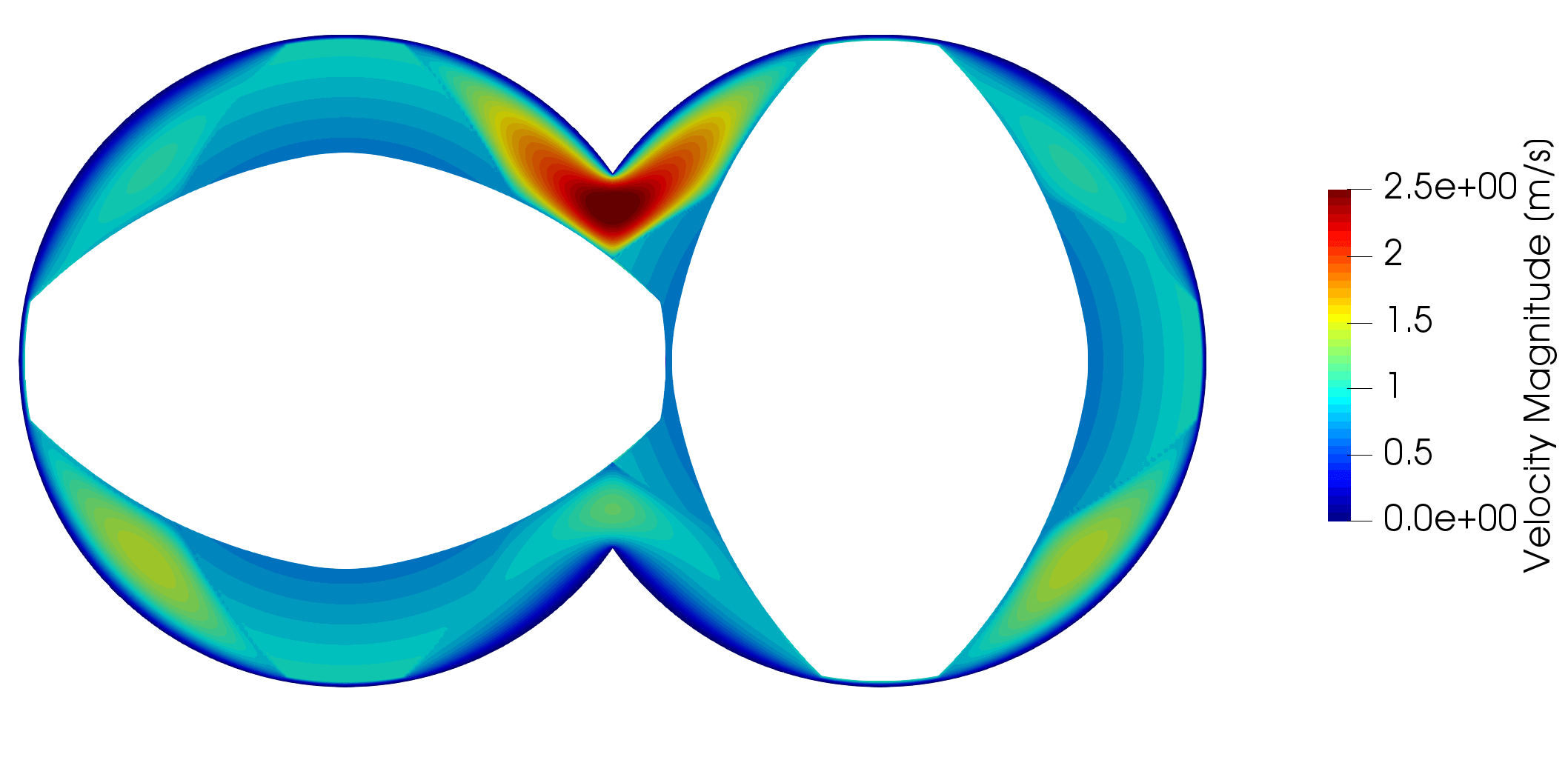}\label{fig:velMag3D2DExtruder2}}
  \caption{Velocity magnitude contour plots on the two sliding interfaces $\Gamma_{SI}^i$ touching $\Gamma_{ff}^{1,2}$ -- (a) front view $0.01mm$ upstream and (b) back view $0.01mm$ downstream of $\Gamma_{ff}^{1,2}$. }
  \label{fig:velMag3D2DExtruder}
\end{figure}

The two individual kneading elements are separately discretized using the Snapping Reference Mesh Update Method (SRMUM) as presented in \cite{helmig2019boundary}.
Thus, we obtain two non-matching surface meshes at the interface between the two individual kneading elements that will be coupled using the presented sliding mesh approach.
The computational domain is the same as used for the motivation in Section \ref{sec:slidingMesh}, see Fig. \ref{fig:slidingMeshExtruder3D}. The resulting interface is shown in Fig. \ref{fig:slidingMeshExtruder2D}. We use the Carreau model in order to account for the non-Newtonian behavior of the plastic melt inside the extruder. The Carreau parameters are given in Table \ref{table:carreau2D}.
We compute the steady flow solution for a screw rotational speed of $\omega_s = 60$ rpm and a pressure difference between inlet and outlet of $\Delta p = 0.5 MPa$. A no-slip condition is used on the barrel, and we set the rotational speed as Dirichlet condition on the screw surface. \\
For the spatial convergence study, we use five different meshes. The mesh parameters are given in Table \ref{table:mesh3DConvergence}: $n_s, n_r$, and $n_a$ describe the number of elements on the screw, in radial, and in axial direction, respectively \cite{helmig2019boundary}. \\
For this test case, we also have to include the weak imposition of boundary conditions at the interface.
The movement of the underlying screw from the other kneading element has to be accurately described on both sides of the interface $\Gamma_{SI} \backslash \Gamma_{ff}$.

\begin{figure}
  \centering
  \subfigure[]{\includegraphics[width=.45\linewidth]{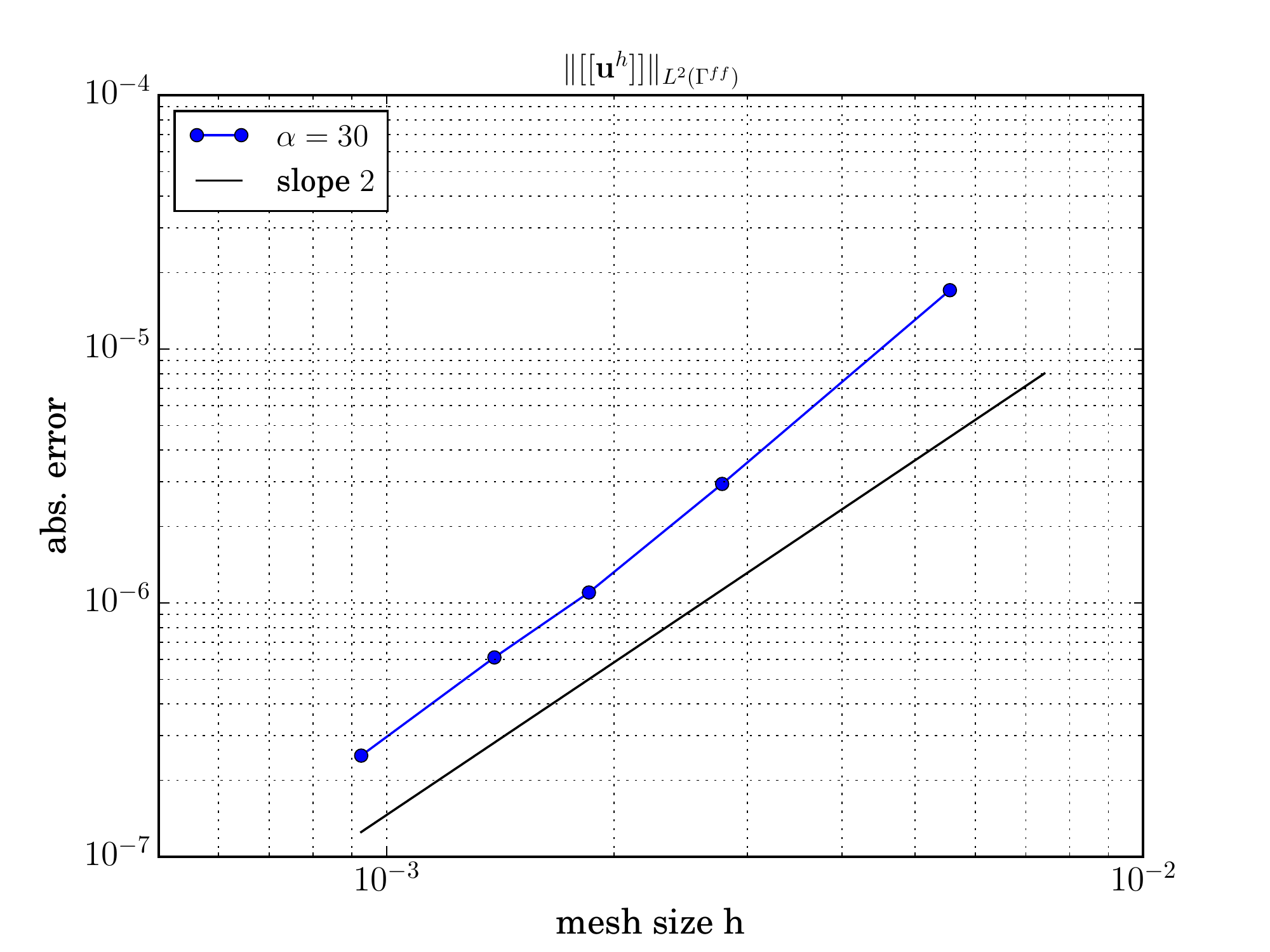}\label{fig:convergence3DExtruderVelocity}}
  \centering
  \subfigure[]{\includegraphics[width=.45\linewidth]{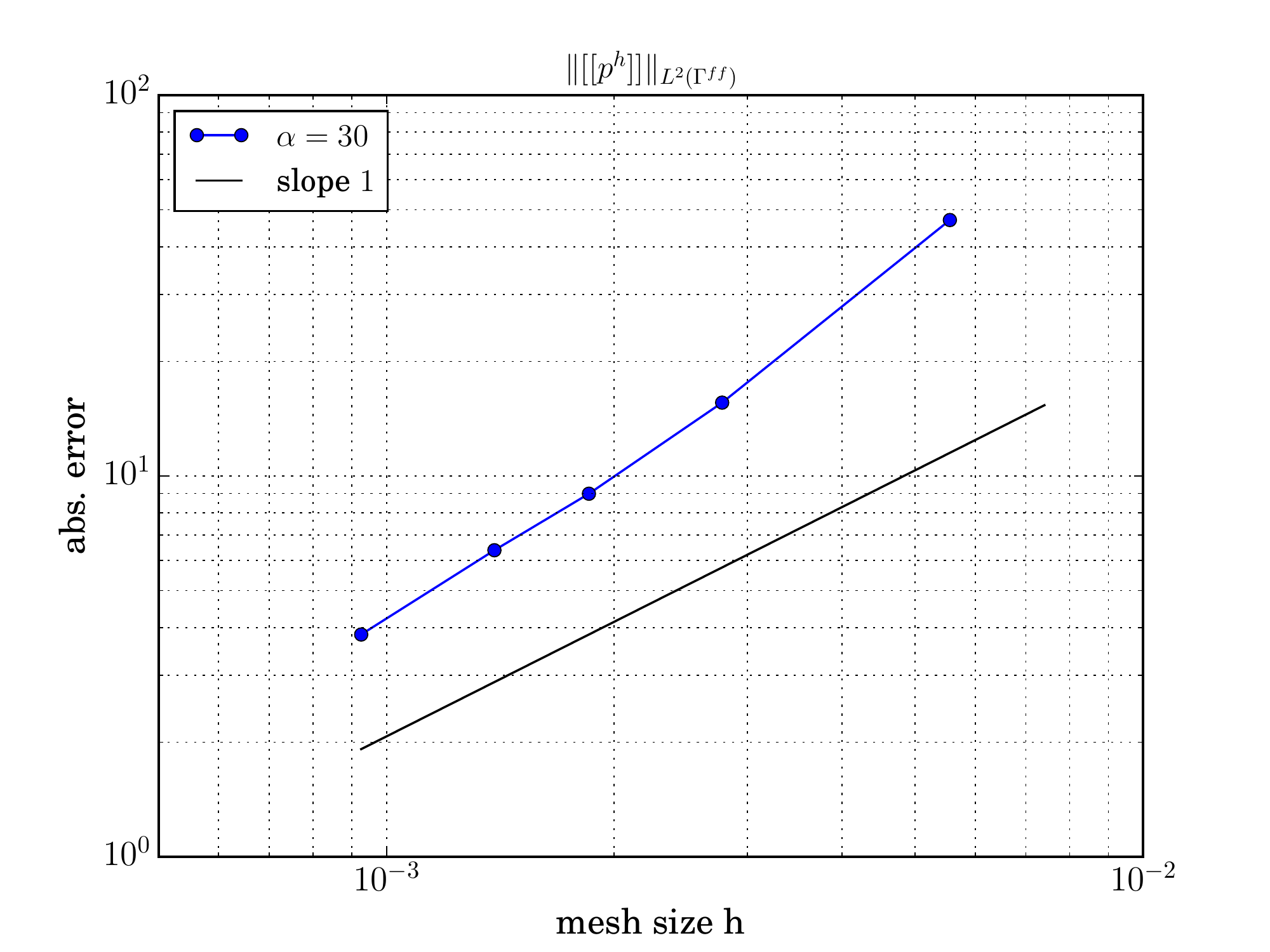}\label{fig:convergence3DExtruderPressure}}
  \caption{Spatial convergence for the flow of plastic melt through a kneading element with two discs staggered by $90^{\circ}$: interface error norms.}
  \label{fig:convergence3DExtruder}
\end{figure}

The solution for the velocity component in $z$-direction on two planes computed on mesh 3 is shown in Fig. \ref{fig:extruder3DzVel}.
We observe that the velocity contours are smooth across the non-matching interface.
Furthermore, we can compare the velocity magnitude on the two sliding interfaces $\Gamma_{SI}^i$, see Fig. \ref{fig:velMag3D2DExtruder}.
We would like to show that the velocity of the opposite screw movement is set correctly.
Looking at the velocity magnitude contours, we observe circular contour lines in those areas that are aligned with the opposite screw body.
This is in perfect agreement with the fact that screw velocity magnitude increases linearly in radial direction of the screw center.
Thus, we can conclude that the weak imposition of Dirichlet boundary condition using the multimesh technique works correctly.

Finally, we want to have a closer look at the interface error norm, which we compute based on Eq. \eqref{eq:interfaceErrorNorm}.
The results for the spatial convergence analysis conducted on five consecutively refined meshes (see Table \ref{table:mesh3DConvergence}) are shown in Fig. \ref{fig:convergence3DExtruder}.
Similar to the results shown in Section \ref{sec:spatialConvergence2D}, we obtain optimal convergence rates for pressure and velocity.
For this sliding mesh setup in 3D using a generalized Newtonian fluid model, this validates the weighting proposed in Section \ref{sec:nitscheCoupling} as well as the scaling for the Nitsche stabilization parameters. This test case shows very good results, so that we can use the proposed method for more complex application cases.

\section{Application Cases} \label{sec:applicationcases}

In the following, we will apply the sliding mesh approach to two relevant physical applications in the plastics manufacturing industry.
The first one is the computation of the temperature-dependent flow of a plastic melt inside a twin-screw extruder with several kneading blocks.
The second one considers the flow inside single-screw extruders with varying design.

\subsection{Temperature-Dependent Flow of Plastic Melt in Twin-Screw Extruder}

We consider the temperature-dependent flow of plastic melt through a twin-screw extruder section with different screw elements.
We will simulate the temperature evolution inside the extruder over several revolutions starting from a constant initial condition.
The non-Newtonian behavior of the plastic melt is modeled by the Cross model with WLF correction, which allows to take the temperature effects into account.
The model parameters are chosen based on a polypropylene from a portfolio of a raw material manufacturer, see Table \ref{table:CrossWLF}.
The melt has density $\rho = 710 \; kg / m^3$, specific heat $c_p = 2900 \; J/kg \;K$, and thermal conductivity $\kappa_0 = 0.2 W/m\;s$.
Similar to \cite{helmig2019boundary}, we use the Prandtl number $Pr = \eta_\infty c_p / \kappa_0$ to relate the momentum diffusivity to thermal diffusivity \cite{sato2004heat}; the Prandtl number used is $Pr =  145 000$.

\begin{figure}[h]
  \centering
  \includegraphics[width=.7\linewidth]{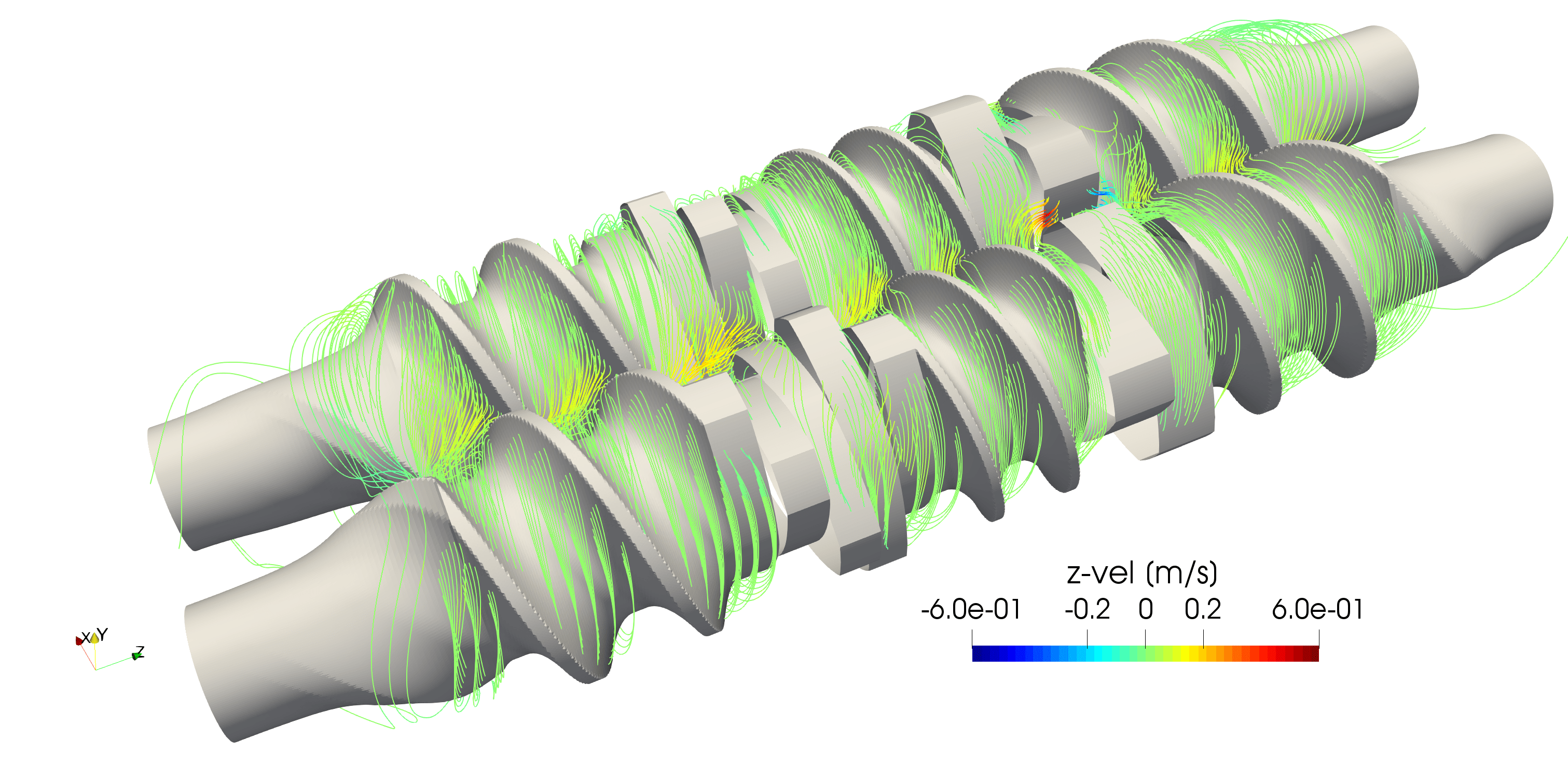}
  \caption{Streamlines for the flow of plastic melt inside a twin-screw extruder with different conveying and kneading elements.}
  \label{fig:extruder3DFullzVel}
\end{figure}

\begin{table}[h]
    \begin{minipage}[t]{0.48\linewidth}
        \centering
        \begin{tabular}{l r}
         \hline
         screw radius $R_s$  & 15.7 $mm$ \\
         center line distance $C_l$ & 26.2 $mm$ \\
         screw-screw clearance $\delta _s$ & 0.3 $mm$ \\
         screw-barrel clearance $\delta _b$ & 0.3 $mm$ \\
         %Pitch length $p_l$ & 28 $mm$ \\
         \hline
         \end{tabular}
         \caption{2D screw geometry parameters.}
         \label{table:screw3DTempUnsteady}
     \end{minipage}
     \begin{minipage}[t]{0.48\linewidth}
         \centering
         \begin{tabular}{l c c}
             \hline
             $D1$  & 1.21e+14 & $Pa \; s$ \\
             $\tau^{*}$ & 25680.0 & $Pa$ \\
             $n$ & 0.2923 & - \\
             $T_{ref}$ & 117.0 & $K$ \\
             $A1$ & 28.32 & - \\
             $A2$ & 51.6 & $K$ \\
             \hline
         \end{tabular}
         \caption{Cross-WLF parameters.}
         \label{table:CrossWLF}
     \end{minipage}
\end{table}

The twin-screw extruder is built up by seven individual screw sections -- artificial relaxation sections at the beginning and the end with a length of $28 \; mm$, where the screw shape is decreased quadratically to a circular shape, two forward-conveying elements with pitch length $p_l = 28 \;mm$, one conveying element with $p_l = 20 \;mm$, one kneading element 45/5/25 (45$^\circ$ staggering angle, 5 discs and $25 \; mm$ length) and one kneading element 60/2/16. The screw centers are at $x=\pm 13.1 \; mm$ and $y=0.0 \; mm$, respectively and the inflow plane is at $z=0 \; mm$.
The 2D screw shape is again based on Booys' law, see Table \ref{table:screw3DTempUnsteady}, and the overall screw setup is shown in Fig. \ref{fig:extruder3DFullzVel}.
All individual screw elements including each disc of a kneading element are discretized using SRMUM.
This allows to update the mesh at each time step in an efficient way without any need for re-meshing.
The number of elements in screw direction and radial direction are $n_s =216$ and $n_r=10$.
The overall number of elements in axial direction is 634.
This has proven to be an appropriate discretization in mesh studies presented in \cite{helmig2019boundary,hinz2019boundary} for similar simulations.
The individual mesh building blocks are coupled at the common interfaces using the presented sliding mesh approach.

\begin{figure}
  \centering
  \subfigure[front view]{\includegraphics[width=.45\linewidth]{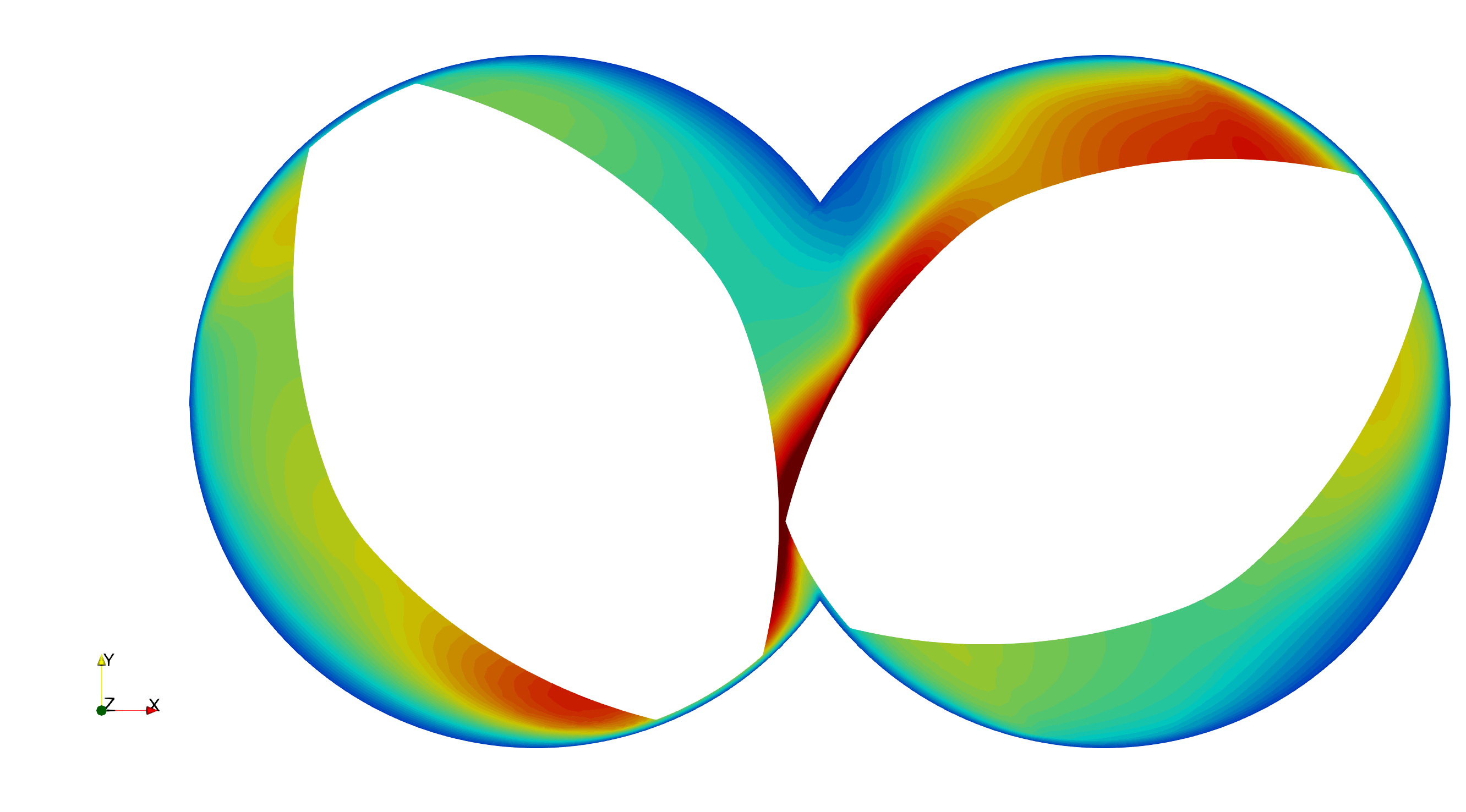}\label{fig:temp3D2DExtruder1}}
  \centering
  \subfigure[back view]{\includegraphics[width=.45\linewidth]{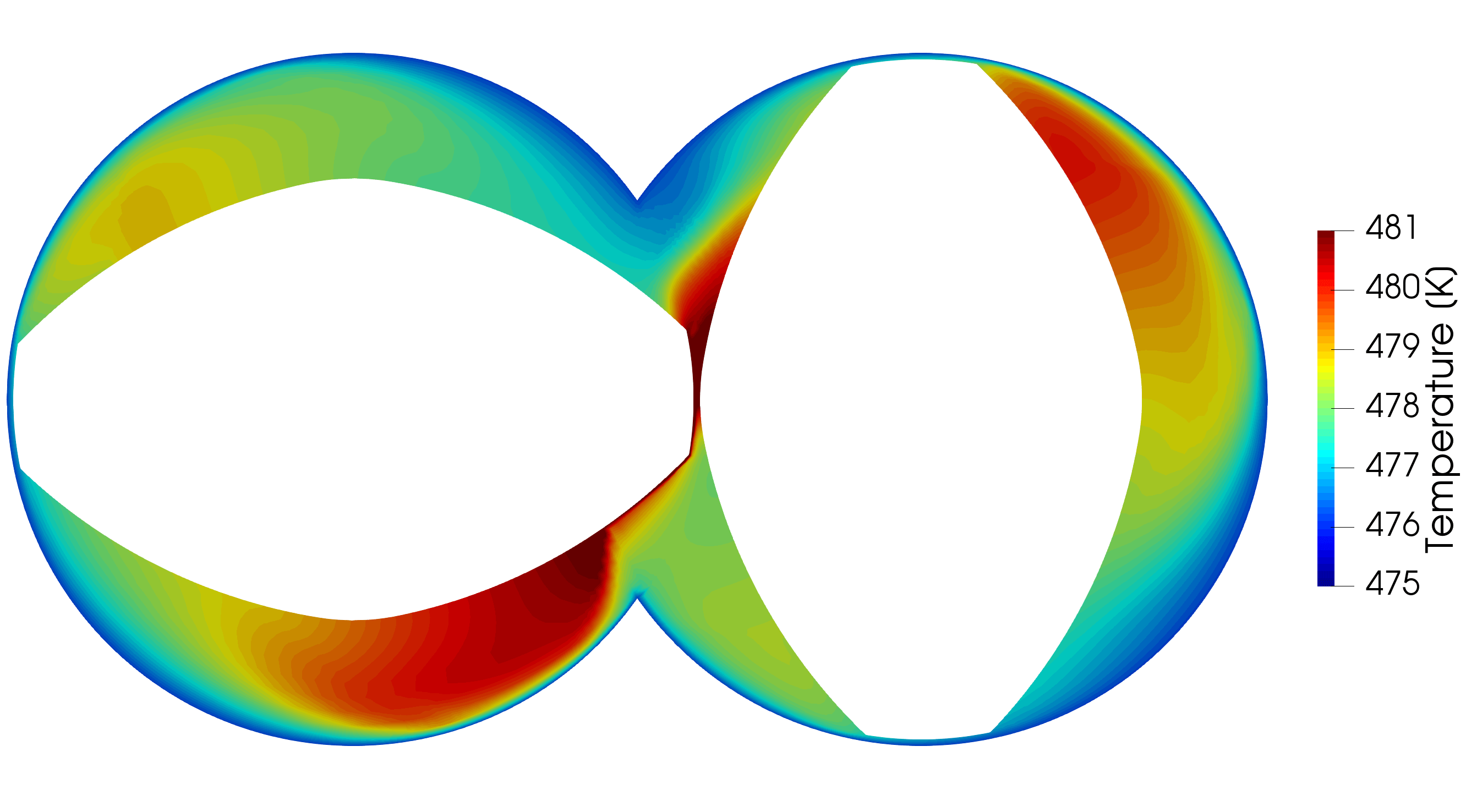}\label{fig:temp3D2DExtruder2}}
  \caption{Temperature contour plots on the two sliding interfaces $\Gamma_{SI}^i$ between the two disc of the kneading element 60/2/16 at $t = 6.54s$ - (a) front view at $z=108.999mm$ and (b) back view at $z=109.001mm$. }
  \label{fig:temp3D2DExtruder}
\end{figure}

\begin{figure}
  \centering
  \subfigure[ $t = 0.62 s$ ]{\includegraphics[width=.45\linewidth]{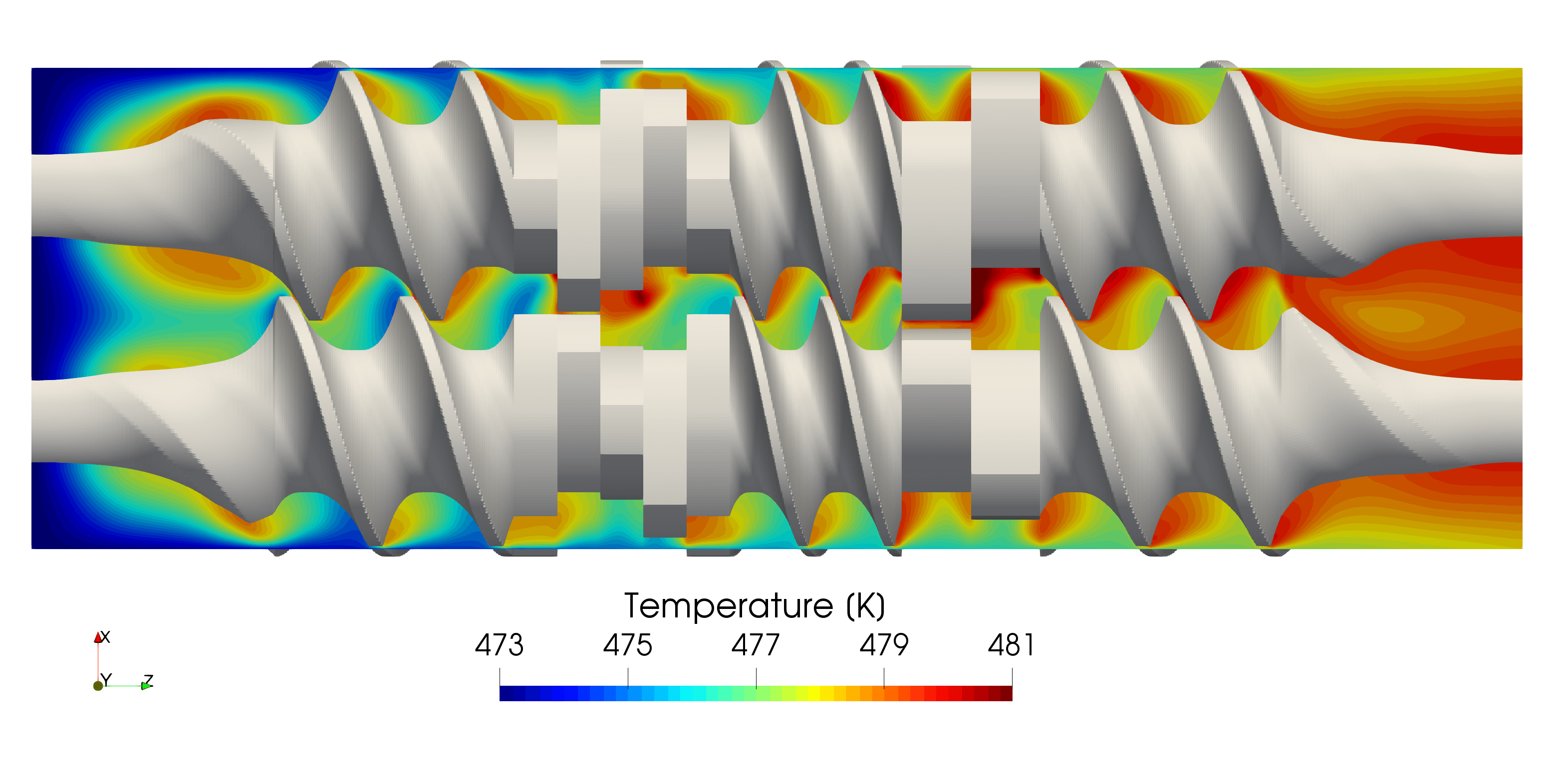}}
  \centering
  \subfigure[$t = 1.0 s$]{\includegraphics[width=.45\linewidth]{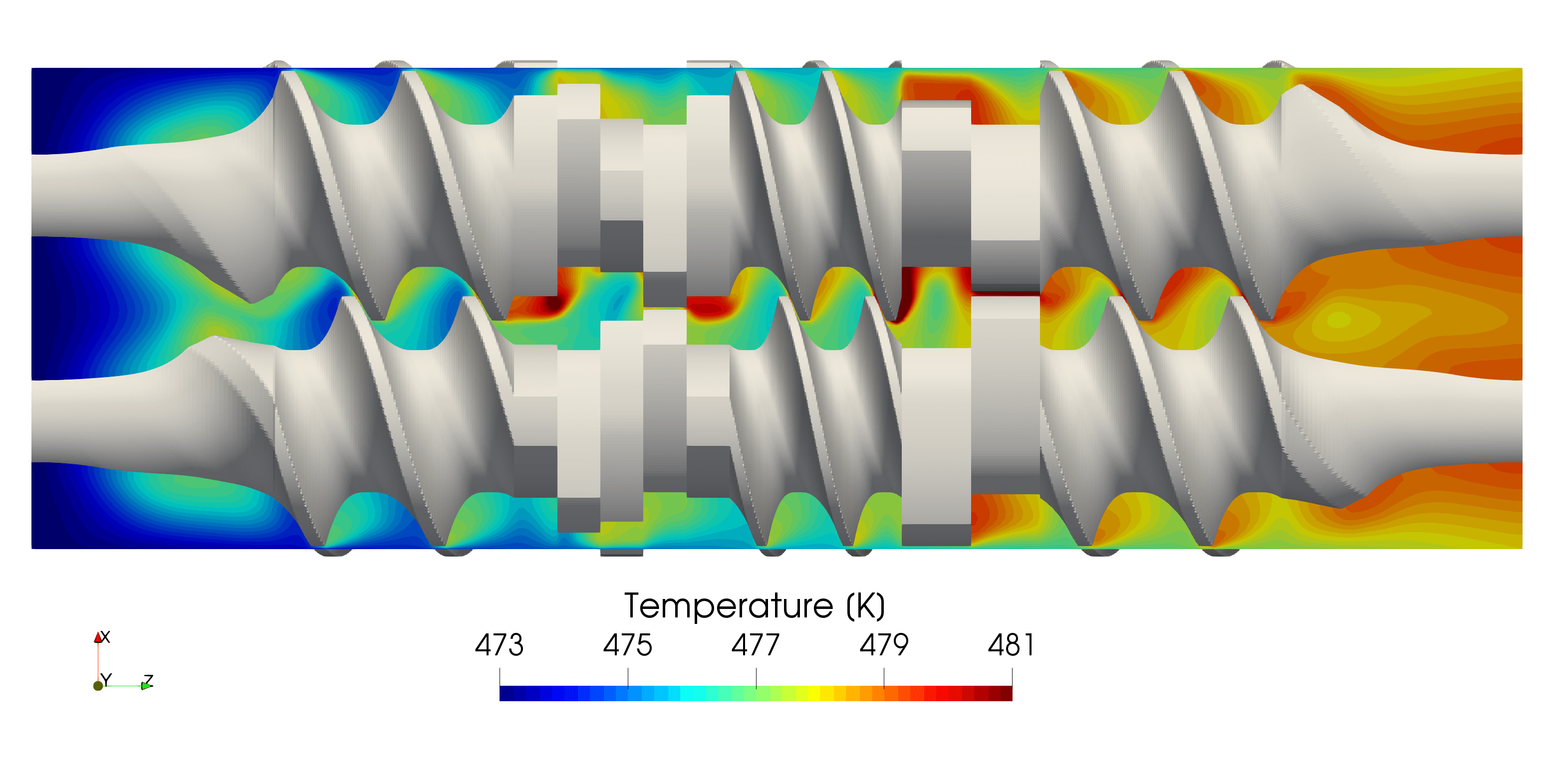}}
  \centering
  \subfigure[$t = 1.62 s$]{\includegraphics[width=.45\linewidth]{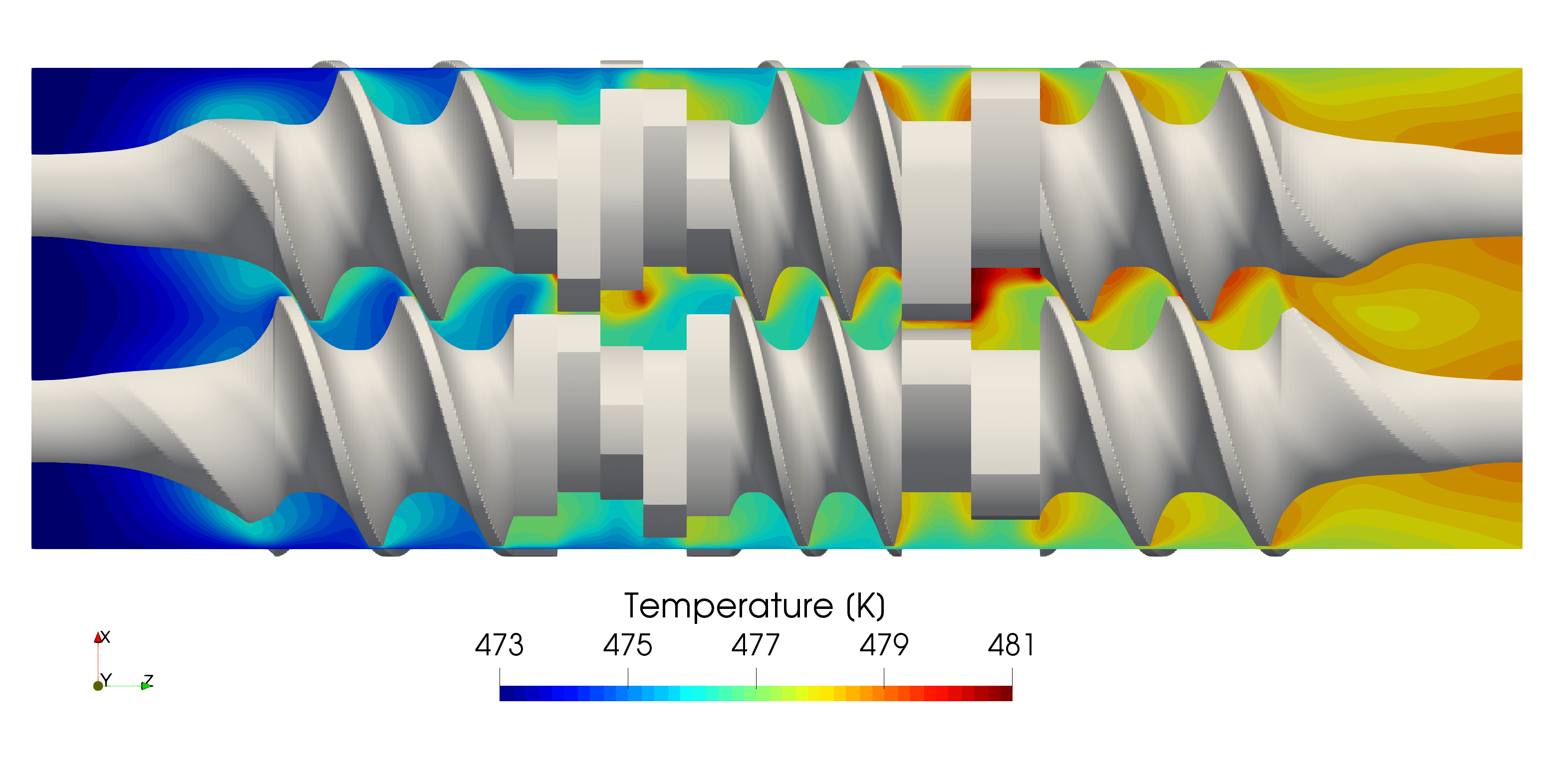}}
  \centering
  \subfigure[$t = 2.0 s$]{\includegraphics[width=.45\linewidth]{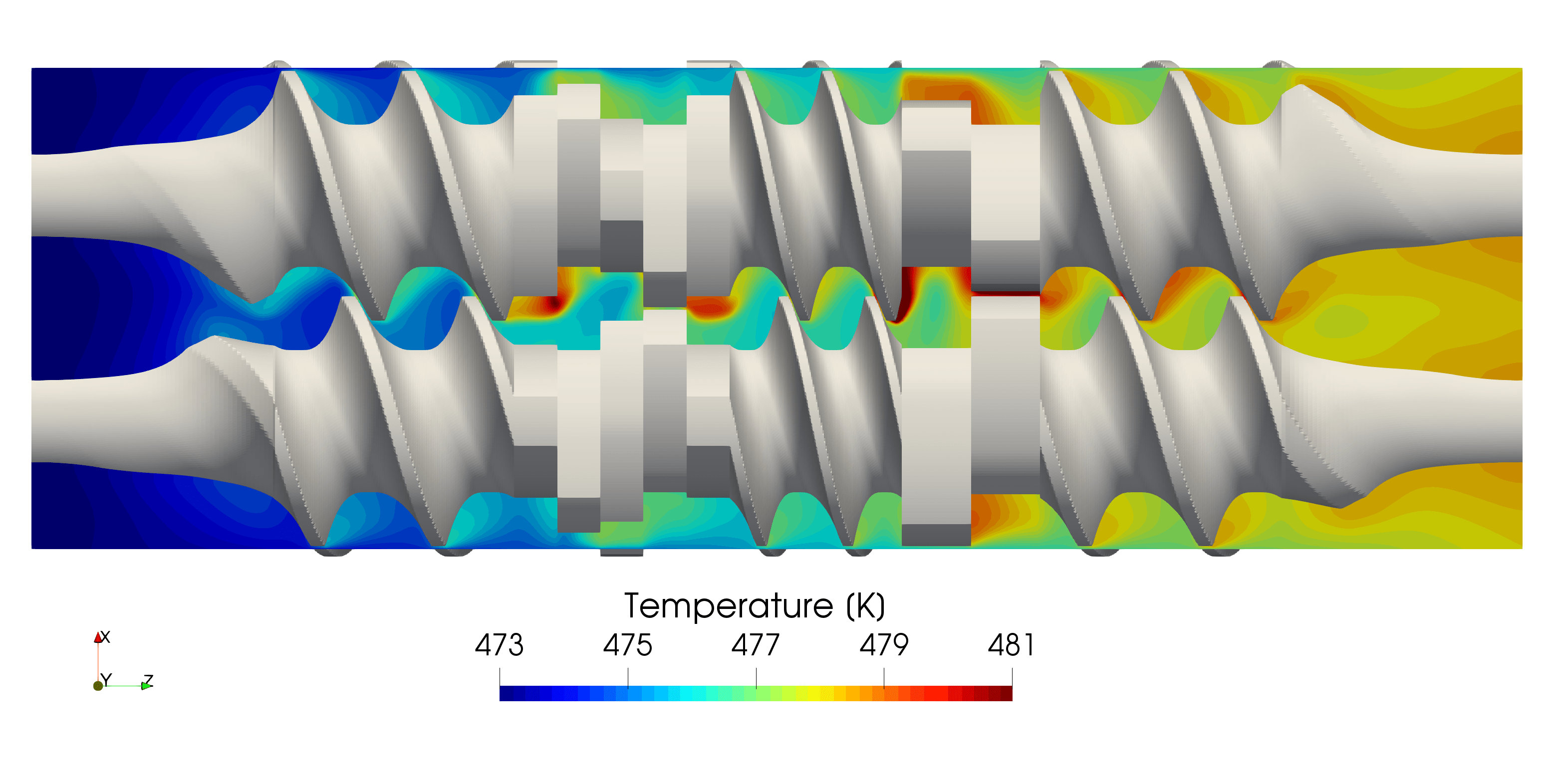}}
  \centering
  \subfigure[$t = 5.62 s$]{\includegraphics[width=.45\linewidth]{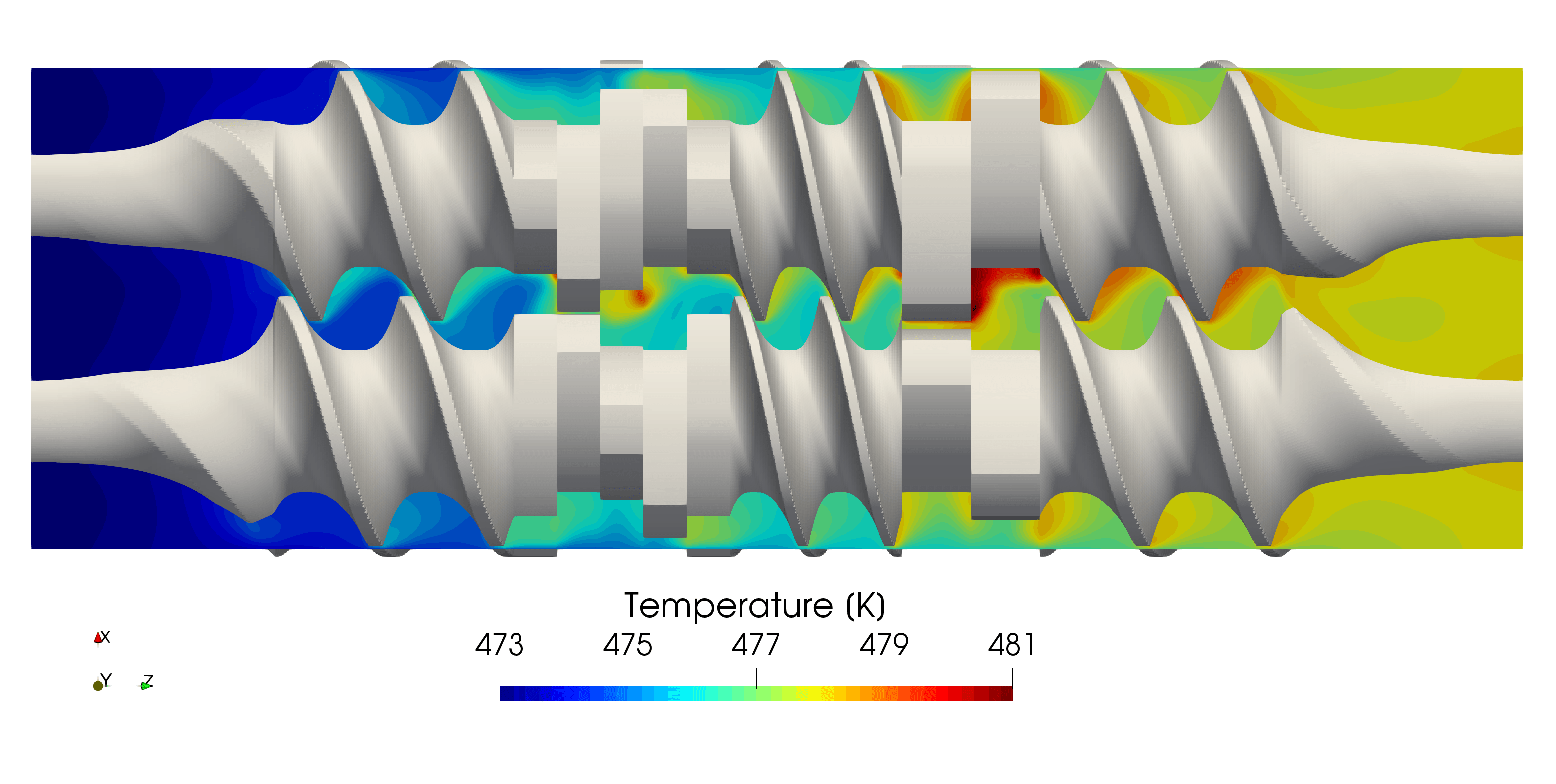}}
  \centering
  \subfigure[$t = 6.0 s$]{\includegraphics[width=.45\linewidth]{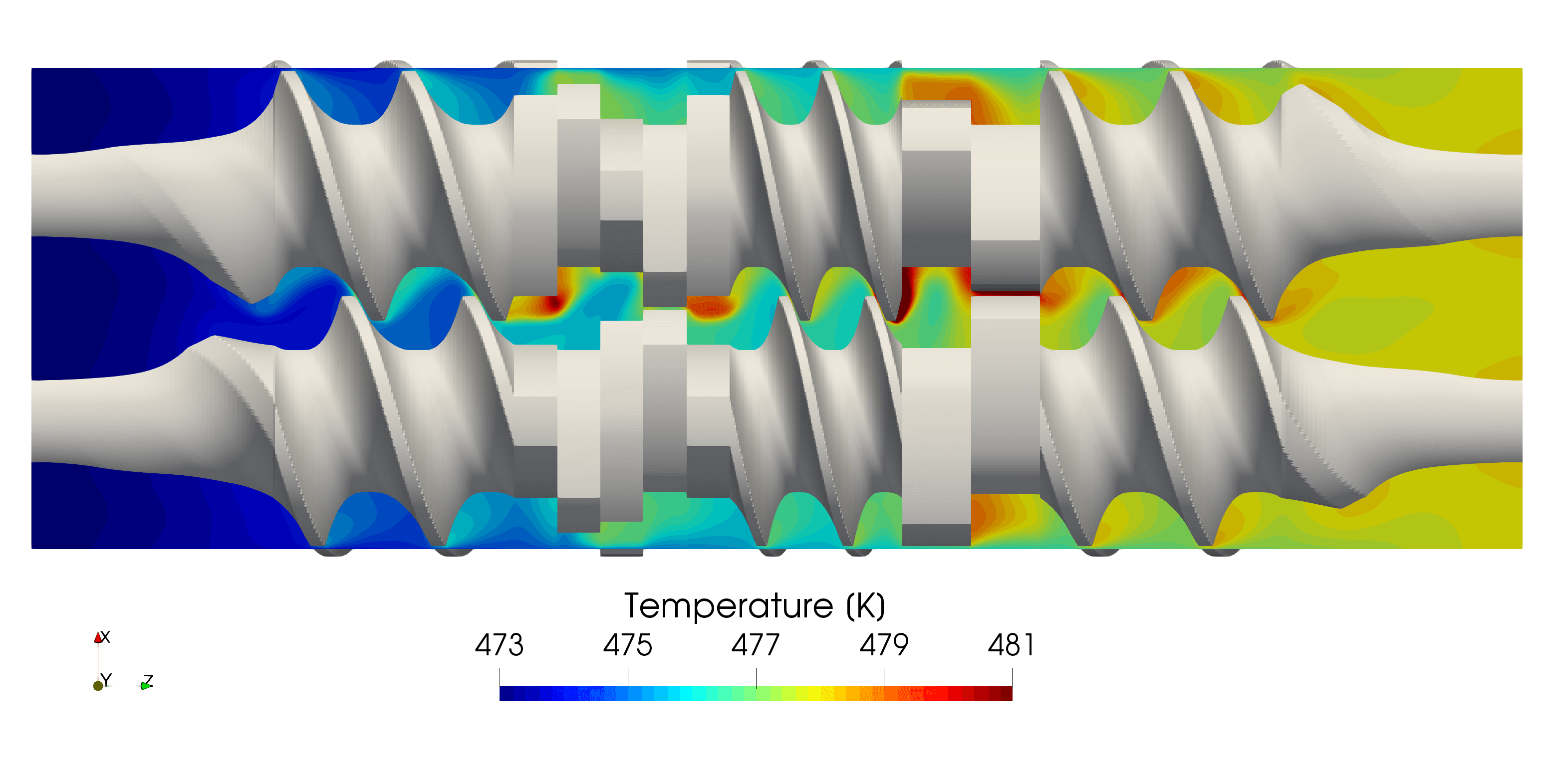}}
  \centering
  \caption{Temperature distribution in the $xz$-plane at $y=6 mm$ inside the twin-screw extruder with different screw sections.}
  \label{fig:temperatureDistributionExtruder3D}
\end{figure}

We simulate a mass flow rate of $\dot{m}= 23.88 kg/h$ which is achieved by setting a uniform Dirichlet inflow condition. A no-slip condition is set at the barrel and a natural boundary condition at the outflow. The screws rotate in mathematically positive direction with $\omega_s = 120$ rpm. The resulting streamlines at $t=6s$ are shown in Fig. \ref{fig:extruder3DFullzVel}. For the temperature, we set a uniform inflow temperate at $T_{in} = 473 K$.
The screws are considered adiabatic and the barrel is heated using a Dirichlet boundary condition with $T_{barrel} = 473 K + 5/z_{max} \cdot 5 K$.
Furthermore, we set the initial temperature of the plastic melt to $T_0 = 480 K$.
The time step size is set to $\Delta t = 0.005s$ based on numerical studies presented in \cite{helmig2019boundary,hinz2019boundary}. The temperature distribution inside the extruder reaches a quasi-steady periodic state after 10 revolutions. \\
Fig. \ref{fig:temp3D2DExtruder} shows the temperature distribution between the two discs of the kneading element 60/2/16 at $t=6.54 s$.
We can observe higher temperatures especially inside the small gap regions due to viscous dissipation.
Furthermore, at the common interface the temperature distributions match well.
This indicates the validity of the presented sliding mesh approach also for the heat equation. \\

Fig. \ref{fig:temperatureDistributionExtruder3D} shows the temperature distribution over time in an $xz$-plane at $y=6mm$.
At $t=0.62s$ and $t=1.0s$ the influence of the initial condition is still visible.
However, cold melt is already pushed into the extruder and a reduction of the temperature due to barrel heating and cooling can be observed.
Additionally, a temperature increase in the vicinity of the small gaps as well as kneading elements due to viscous heating is visible.
The overall temperature decreases over time compared to the initial condition. The temperature distribution after $t=1s$ is already close to the quasi-steady state.
The biggest differences occur in the outflow region. It is noteworthy, that the average temperature in any cross section is higher than the given temperature at the corresponding barrel position.
This clearly shows the importance of viscous heating effects in the extruder and is also illustrated in Fig. \ref{fig:temp3D2DExtruder}.
The computed temperature distributions are also qualitatively in accordance to steady-state results for different kneading elements presented in \cite{kalyon2007integrated}. \\
Using the sliding mesh for this complex application shows its potential to couple individual moving boundary-conforming meshes at common interfaces.

\subsection{Single-Screw Extruder}

\begin{figure}[h]
  \centering
  \includegraphics[width=.8\linewidth]{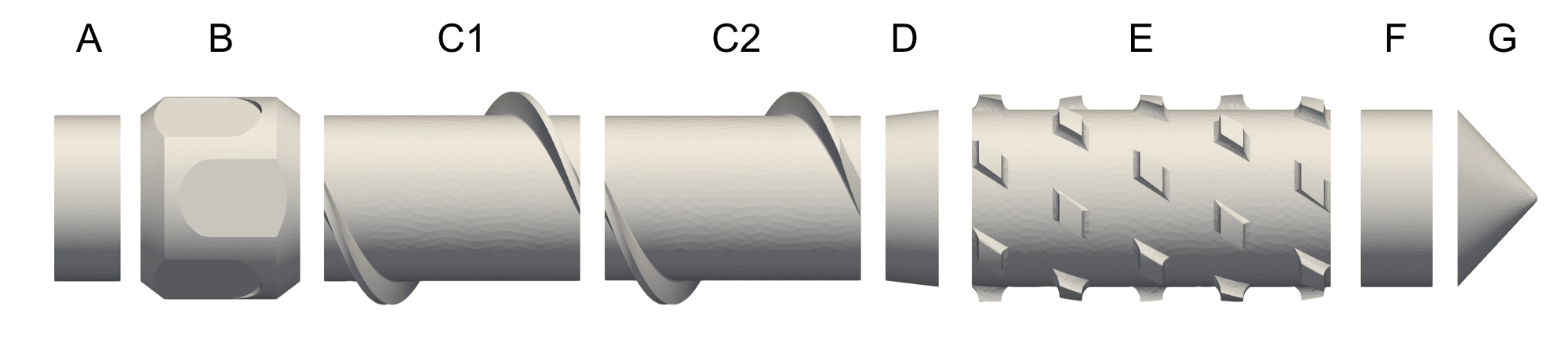}
  \caption{Different individual screw parts for a single-screw extruder.}
  \label{fig:SSEParts}
\end{figure}

In the following, we will apply the sliding mesh approach to a different application -- still in extruding -- but this time we consider single screw extruders (SSE).
Similar to twin-screw extruders, single-screw extruders are built up by combining different screw geometries that each serve a different purpose. A selection of potential screw parts is shown in Fig. \ref{fig:SSEParts}.
Conveying elements (C) are used in order to transport the melt forward, whereas Maddock elements (B) serve as barriers that disperse the melt and decrease the pressure. In contrast to twin-screw extruders, the conveying elements do not provide enough mixing. Thus, extra distributive mixing elements (E) have to be used. Finally, there is a section (G) that leads to the extrusion die, and other transition sections.
An open research question is how to assemble the individual screw shapes in an optimal order, considering for example the mixing behavior or pressure loss, to name two potential design objectives \cite{eusterholz2018cfd}. The screw designs again feature very small gaps between screw and barrel, e.g., with a barrel diameter of $60\;mm$ the smallest gap size between the conveying screw element and the barrel is only $0.3\; mm$.
Thus, it is again of utmost importance to have a good discretization in those regions in order to capture all flow effects.
This is where the sliding mesh approach comes into play.
It allows to create individual boundary-conforming meshes for each individual screw section and to couple them at the common interface.
Thus, it is not necessary to re-mesh in case one assembles the screws in a different way.
In the following, we aim to demonstrate the usability of the sliding mesh approach to compute the flow field for four different screw combinations.

\begin{table}[h]
    \centering
    \begin{tabular}{l c c c c}
        \hline
        & 1 & 2 & 3 & 4 \\
      combinations & A-C1-B-D-E-F-G & A-C1-C2-B-D-E-F-G & A-B-C1-C2-D-E-F-G & A-B-C1-D-E-F-G \\
      elements & 6,847,102 & 10,192,464 & 10,192,464 & 6,847,102 \\
      \hline
    \end{tabular}
    \caption{Screw combinations based on Fig. \ref{fig:SSEParts} for four different screw combinations as well as the total number of elements for each resulting mesh.}
    \label{table:SSECombinations}
\end{table}

The screw configurations are given in Table \ref{table:SSECombinations}. The barrel diameter is $60 \;mm$. As fluid, we consider corn syrup that can be modeled as a Newtonian fluid with density $\rho = 1400 kg/m^3$ and viscosity $\eta = 4.7 Pa \; s$. The screw rotates with rotational speed of $\omega_s = 60$ rpm. We set natural boundary conditions at the inlet and outlet. Thus, the flow rate will be mainly determined by the transport capacity of the conveying elements in combination with the Maddock elements. For simplicity, we only simulate a steady-state flow without considering temperature effects.

\begin{figure}[h]
  \centering
  \includegraphics[width=.9\linewidth]{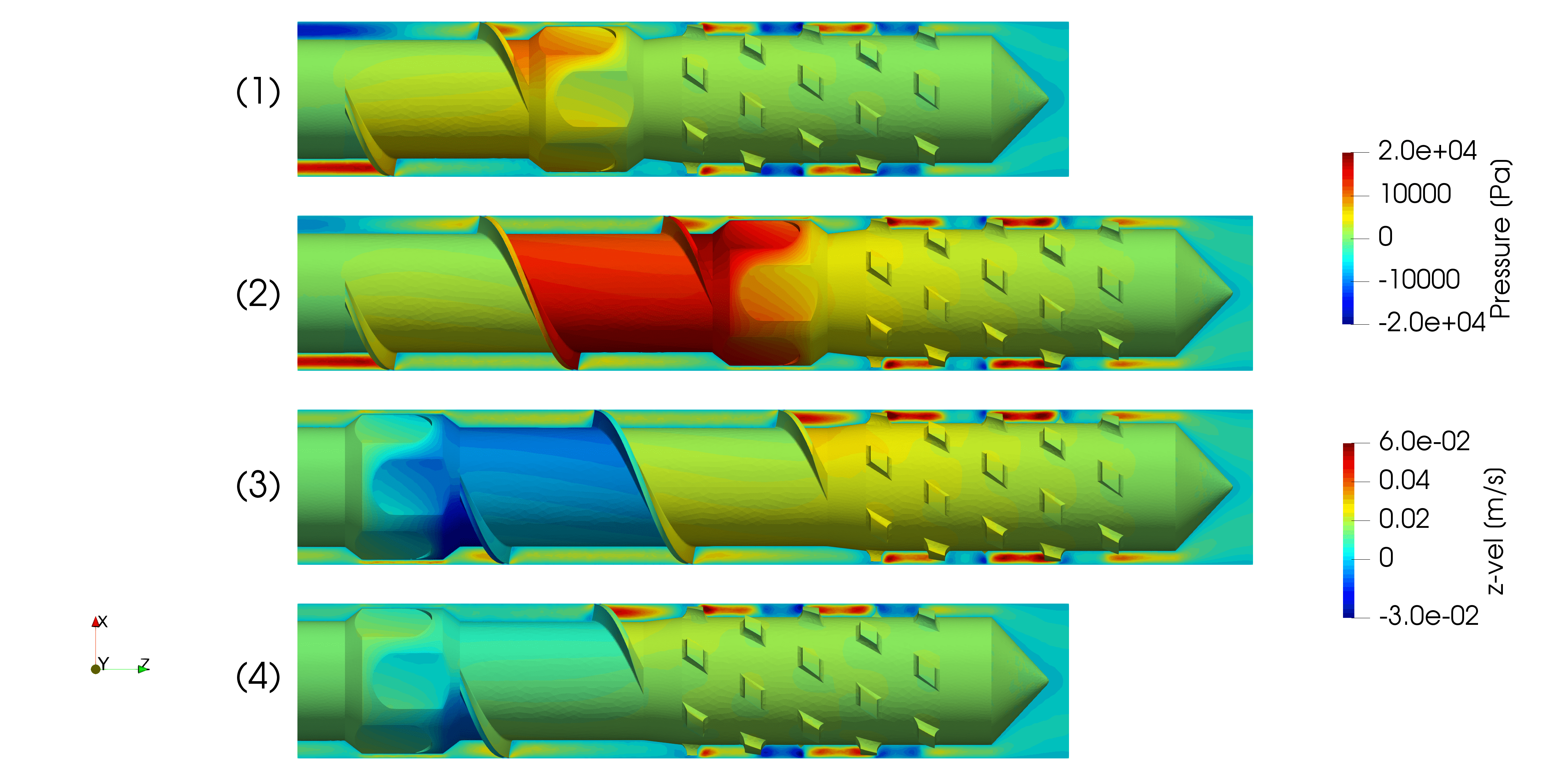}
  \caption{Pressure field on screw surface and $z$-velocity field in $xz$-plane through screw center for four different screw configurations.}
  \label{fig:SSE3DPressurezVel}
\end{figure}

\begin{table}[h!]
    \centering
    \begin{tabular}{l c c  c c}
        \hline
        & 1 & 2 & 3 & 4 \\
      $\dot{m} \; [ kg/h ]$ & $37.34 $ & $70.76 $ & $72.53 $ & $40.25$ \\
      \hline
    \end{tabular}
    \caption{Mass flow rates for 4 screw configurations.}
    \label{table:SSEMassFlowRate}
\end{table}

The resulting mass flow rates are given in Table \ref{table:SSEMassFlowRate}. Comparing configuration 1 and 4, as well as configuration 2 and 3, we can observe that there is a difference in the mass flow depending on the position of the Maddock element and the conveying elements. Adding an additional conveying element nearly doubles the mass flow rate, which is physically expected since the flow rate is mainly driven by the conveying elements.
Furthermore, a factor that is slightly smaller than two makes sense since the mixing element also features a certain conveying behavior.
The absolute increase of the mass flow rate by adding an additional conveying element is also very similar for configuration 2 ($ 33.42 \; kg/h$) and configuration 3 ($32.28 \; kg/h $).
Fig. \ref{fig:SSE3DPressurezVel} shows the flow for all four configuration.
The pressure field is shown on the screw surface and the velocity in $z$-direction on the plane $xz$-plane through the screw axis. The velocity as well as pressure contours are smooth across the non-matching interfaces. Furthermore, the mixing effects inside the mixing elements are clearly discernable by high negative velocities in $z$-direction. Additionally, the velocity magnitude reflects the different flow rates in case of one or two conveying elements.
The pressure contours also demonstrate the pressure decrease over the Maddock element, which is mainly compensated by the conveying elements. \\

% !TEX root = ./main.tex
%!TeX spellcheck = en-US
\section{Conclusion} \label{sec:conclusion}

Within this work, we presented a stabilized finite element formulation for fluid and temperature equations that allows to couple sliding domains with boundary-conforming meshes at common interfaces.
This approach allows to use highly-optimized boundary-conforming meshes for individual subdomains that can then be coupled weakly at their -- possible only partially overlapping -- common interface.
We solve the non-isothermal, incompressible Navier-Stokes equations on the individual subdomains for Newtonian as well as for generalized Newtonain fluid models.
The formulation is stabilized using residual based stabilization technique.
The solution fields of each individual subdomain are coupled weakly.
The interface coupling conditions are enforced using Nitsche's method.
Additionally, the method imposes Dirichlet conditions weakly on boundaries that partly slide over rigid segments. \\
We verified the method for the 2D Taylor-Green flow for viscous as well as convective flow regimes.
Optimal convergence rates were obtained in a spatial convergence study, which suggests the correct scaling of our Nitsche stabilization terms.
In a second validation step, a 3D test case considering the flow of plastic melt inside a twin-screw extruder kneading element with two discs was used to validate the method with respect to: (1) 3D cases, (2) generalized Newtonian models, and (3) the weak imposition of Dirichlet boundary conditions on the sliding interfaces. Again, optimal convergence rates for the interface coupling could be obtained in a spatial convergence study. \\
First steps towards relevant industrial applications have been also made.
We computed the time-dependent temperature and flow field of plastic melt inside a twin-screw extruder with several kneading and conveying elements.
The SRMUM was used as mesh update method that allows to use boundary-conforming meshes for the individual rotating screw elements. The boundary conforming meshes were only coupled at the non-matching interface between individual kneading discs.
As a second industrially relevant example, we computed the flow inside a single-screw extruder with different screw sections.
Individual screw sections were meshed independently. We assembled the sections in four different ways. The resulting flow fields showed the expected behavior. This test case demonstrated the applicability of the presented approach to discrete optimization of single-screw extruders without any need for re-meshing. \\
We believe that the presented method makes use of the benefits of both the boundary-conforming and unfitted approaches in a way tailored to the considered application cases. In the future, the method will be applied to more complex application cases, including realistic design optimization of extruders.

\section*{Acknowledgements}

The computations were conducted on computing clusters supplied by the J\"ulich Aachen Research Alliance (JARA) and the IT Center of the RWTH Aachen University. Furthermore, the authors would like to thank Andre Massing for valuable input concerning Nitsche's method in outside coffee breaks during the rare sunny summer days in Aachen.

\bibliographystyle{elsarticle-num}
\bibliography{bibliography}

\begin{thebibliography}{10}
\expandafter\ifx\csname url\endcsname\relax
  \def\url#1{\texttt{#1}}\fi
\expandafter\ifx\csname urlprefix\endcsname\relax\def\urlprefix{URL }\fi
\expandafter\ifx\csname href\endcsname\relax
  \def\href#1#2{#2} \def\path#1{#1}\fi

\bibitem{de2007mesh}
A.~De~Boer, M.~Van~der Schoot, H.~Bijl, Mesh deformation based on radial basis
  function interpolation, Computers \& structures 85~(11-14) (2007) 784--795.

\bibitem{rendall2009efficient}
T.~C. Rendall, C.~B. Allen, Efficient mesh motion using radial basis functions
  with data reduction algorithms, Journal of Computational Physics 228~(17)
  (2009) 6231--6249.

\bibitem{johnson1994mesh}
A.~A. Johnson, T.~E. Tezduyar, Mesh update strategies in parallel finite
  element computations of flow problems with moving boundaries and interfaces,
  Computer methods in applied mechanics and engineering 119~(1-2) (1994)
  73--94.

\bibitem{stein2003mesh}
K.~Stein, T.~Tezduyar, R.~Benney, Mesh moving techniques for fluid-structure
  interactions with large displacements, J. Appl. Mech. 70~(1) (2003) 58--63.

\bibitem{takizawa2020low}
K.~Takizawa, T.~E. Tezduyar, R.~Avsar, A low-distortion mesh moving method
  based on fiber-reinforced hyperelasticity and optimized zero-stress state,
  Computational Mechanics (2020) 1--25.

\bibitem{behr1999shear}
M.~Behr, T.~Tezduyar, The shear-slip mesh update method, Computer Methods in
  Applied Mechanics and Engineering 174~(3-4) (1999) 261--274.

\bibitem{behr2003shear}
M.~Behr, D.~Arora, Shear-slip mesh update method: Implementation and
  applications, Computer Methods in Biomechanics and Biomedical Engineering
  6~(2) (2003) 113--123.

\bibitem{key2018virtual}
F.~Key, L.~Pauli, S.~Elgeti, The virtual ring shear-slip mesh update method,
  Computers \& Fluids 172 (2018) 352--361.

\bibitem{helmig2019boundary}
J.~Helmig, M.~Behr, S.~Elgeti, Boundary-conforming finite element methods for
  twin-screw extruders: Unsteady-temperature-dependent-non-newtonian
  simulations, Computers \& Fluids 190 (2019) 322--336.

\bibitem{hinz2019boundary}
J.~Hinz, J.~Helmig, M.~M{\"o}ller, S.~Elgeti, Boundary-conforming finite
  element methods for twin-screw extruders using spline-based parameterization
  techniques, Computer Methods in Applied Mechanics and Engineering (2019)
  112740.

\bibitem{rane2013grid}
S.~Rane, A.~Kovacevic, N.~Stosic, M.~Kethidi, Grid deformation strategies for
  cfd analysis of screw compressors, International Journal of Refrigeration
  36~(7) (2013) 1883--1893.

\bibitem{lehrenfeld2017higher}
C.~Lehrenfeld, A higher order isoparametric fictitious domain method for level
  set domains, in: Geometrically Unfitted Finite Element Methods and
  Applications, Springer, 2017, pp. 65--92.

\bibitem{peskin2002immersed}
C.~S. Peskin, The immersed boundary method, Acta numerica 11 (2002) 479--517.

\bibitem{steger1991chimera}
J.~Steger, The chimera method of flow simulation, Vol. 188, Workshop on applied
  CFD, Univ of Tennessee Space Institute, 1991.

\bibitem{belk1995role}
D.~M. Belk, The role of overset grids in the development of the general purpose
  cfd code, in: NASA CONFERENCE PUBLICATION, NASA, 1995, pp. 193--193.

\bibitem{houzeaux2003chimera}
G.~Houzeaux, R.~Codina, A chimera method based on a dirichlet/neumann (robin)
  coupling for the navier--stokes equations, Computer Methods in Applied
  Mechanics and Engineering 192~(31-32) (2003) 3343--3377.

\bibitem{houzeaux2017domain}
G.~Houzeaux, J.~Cajas, M.~Discacciati, B.~Eguzkitza, A.~Gargallo-Peir{\'o},
  M.~Rivero, M.~V{\'a}zquez, Domain decomposition methods for domain
  composition purpose: Chimera, overset, gluing and sliding mesh methods,
  Archives of Computational Methods in Engineering 24~(4) (2017) 1033--1070.

\bibitem{gerstenberger2008extended}
A.~Gerstenberger, W.~A. Wall, An extended finite element method/lagrange
  multiplier based approach for fluid--structure interaction, Computer Methods
  in Applied Mechanics and Engineering 197~(19-20) (2008) 1699--1714.

\bibitem{mayer20103d}
U.~M. Mayer, A.~Popp, A.~Gerstenberger, W.~A. Wall, 3d fluid--structure-contact
  interaction based on a combined xfem fsi and dual mortar contact approach,
  Computational Mechanics 46~(1) (2010) 53--67.

\bibitem{fard2012extended}
A.~S. Fard, M.~A. Hulsen, P.~D. Anderson, Extended finite element method for
  viscous flow inside complex three-dimensional geometries with moving internal
  boundaries, International Journal for Numerical Methods in Fluids 70~(6)
  (2012) 775--792.

\bibitem{villanueva2014density}
C.~H. Villanueva, K.~Maute, Density and level set-xfem schemes for topology
  optimization of 3-d structures, Computational Mechanics 54~(1) (2014)
  133--150.

\bibitem{nitsche1971variationsprinzip}
J.~Nitsche, {\"U}ber ein variationsprinzip zur l{\"o}sung von
  dirichlet-problemen bei verwendung von teilr{\"a}umen, die keinen
  randbedingungen unterworfen sind, in: Abhandlungen aus dem mathematischen
  Seminar der Universit{\"a}t Hamburg, Vol.~36, Springer, 1971, pp. 9--15.

\bibitem{hansbo2002unfitted}
A.~Hansbo, P.~Hansbo, An unfitted finite element method, based on nitsche’s
  method, for elliptic interface problems, Computer methods in applied
  mechanics and engineering 191~(47-48) (2002) 5537--5552.

\bibitem{massing2014stabilized}
A.~Massing, M.~G. Larson, A.~Logg, M.~E. Rognes, A stabilized nitsche
  overlapping mesh method for the stokes problem, Numerische Mathematik 128~(1)
  (2014) 73--101.

\bibitem{schott2016stabilized}
B.~Schott, S.~Shahmiri, R.~Kruse, W.~Wall, A stabilized nitsche-type extended
  embedding mesh approach for 3d low-and high-reynolds-number flows,
  International Journal for Numerical Methods in Fluids 82~(6) (2016) 289--315.

\bibitem{johansson2019multimesh}
A.~Johansson, B.~Kehlet, M.~G. Larson, A.~Logg, Multimesh finite element
  methods: Solving pdes on multiple intersecting meshes, Computer Methods in
  Applied Mechanics and Engineering 343 (2019) 672--689.

\bibitem{dokken2019multimesh}
J.~S. Dokken, A.~Johansson, A.~Massing, S.~W. Funke, A multimesh finite element
  method for the navier-stokes equations based on projection methods, arXiv
  preprint arXiv:1912.06392.

\bibitem{parvizian2007finite}
J.~Parvizian, A.~D{\"u}ster, E.~Rank, Finite cell method, Computational
  Mechanics 41~(1) (2007) 121--133.

\bibitem{hoang2017mixed}
T.~Hoang, C.~V. Verhoosel, F.~Auricchio, E.~H. van Brummelen, A.~Reali, Mixed
  isogeometric finite cell methods for the stokes problem, Computer Methods in
  Applied Mechanics and Engineering 316 (2017) 400--423.

\bibitem{ianus2014mesh}
O.~Mierka, T.~Theis, T.~Herken, S.~Turek, V.~Sch{\"o}ppner, F.~Platte, Mesh
  Deformation Based Finite Element-Fictitious Boundary Method (FEM-FBM) for the
  Simulation of Twin-screw Extruders, Citeseer, 2014.

\bibitem{bazilevs2008nurbs}
Y.~Bazilevs, T.~Hughes, Nurbs-based isogeometric analysis for the computation
  of flows about rotating components, Computational Mechanics 43~(1) (2008)
  143--150.

\bibitem{hsu2012fluid}
M.-C. Hsu, Y.~Bazilevs, Fluid--structure interaction modeling of wind turbines:
  simulating the full machine, Computational Mechanics 50~(6) (2012) 821--833.

\bibitem{bazilevs2007weak}
Y.~Bazilevs, T.~J. Hughes, Weak imposition of dirichlet boundary conditions in
  fluid mechanics, Computers \& Fluids 36~(1) (2007) 12--26.

\bibitem{massing2018stabilized}
A.~Massing, B.~Schott, W.~A. Wall, A stabilized nitsche cut finite element
  method for the oseen problem, Computer Methods in Applied Mechanics and
  Engineering 328 (2018) 262--300.

\bibitem{donea2003finite}
J.~Donea, A.~Huerta, Finite element methods for flow problems, John Wiley \&
  Sons, 2003.

\bibitem{forster2006geometric}
C.~F{\"o}rster, W.~A. Wall, E.~Ramm, On the geometric conservation law in
  transient flow calculations on deforming domains, International Journal for
  Numerical Methods in Fluids 50~(12) (2006) 1369--1379.

\bibitem{carreau1979review}
P.~Carreau, D.~De~Kee, Review of some useful rheological equations, The
  Canadian Journal of Chemical Engineering 57~(1) (1979) 3--15.

\bibitem{rudolph2014polymer}
N.~Rudolph, T.~A. Osswald, Polymer rheology: fundamentals and applications,
  Carl Hanser Verlag GmbH Co KG, 2014.

\bibitem{pauli2017stabilized}
L.~Pauli, M.~Behr, On stabilized space-time fem for anisotropic meshes:
  Incompressible navier--stokes equations and applications to blood flow in
  medical devices, International Journal for Numerical Methods in Fluids 85~(3)
  (2017) 189--209.

\bibitem{hughes2018multiscale}
T.~J. Hughes, G.~Scovazzi, L.~P. Franca, Multiscale and stabilized methods,
  Encyclopedia of Computational Mechanics Second Edition (2018) 1--64.

\bibitem{jansen1999better}
K.~E. Jansen, S.~S. Collis, C.~Whiting, F.~Shakib, A better consistency for
  low-order stabilized finite element methods, Computer methods in applied
  mechanics and engineering 174~(1-2) (1999) 153--170.

\bibitem{pauli2016stabilized}
L.~Pauli, Stabilized Finite Element Methods for Computational Design of
  Blood-Handling Devices, Verlag Dr. Hut, 2016.

\bibitem{cgal:eb-20a}
{The CGAL Project},
  \href{https://doc.cgal.org/5.0.1/Manual/packages.html}{{CGAL} User and
  Reference Manual}, {5.0.1} Edition, {CGAL Editorial Board}, 2020.
\newline\urlprefix\url{https://doc.cgal.org/5.0.1/Manual/packages.html}

\bibitem{pearson1964computational}
C.~E. Pearson, A computational method for time-dependant two-dimensional
  incompressible viscous flow problems, Sperry Rand Research Centre, 1964.

\bibitem{forti2015semi}
D.~Forti, L.~Ded{\`e}, Semi-implicit bdf time discretization of the
  navier--stokes equations with vms-les modeling in a high performance
  computing framework, Computers \& Fluids 117 (2015) 168--182.

\bibitem{burman2014fictitious}
E.~Burman, P.~Hansbo, Fictitious domain methods using cut elements: Iii. a
  stabilized nitsche method for stokes’ problem, ESAIM: Mathematical
  Modelling and Numerical Analysis 48~(3) (2014) 859--874.

\bibitem{sarhangi2012adaptive}
A.~S. Fard, M.~Hulsen, H.~Meijer, N.~Famili, P.~Anderson, Adaptive
  non-conformal mesh refinement and extended finite element method for viscous
  flow inside complex moving geometries, International Journal for Numerical
  Methods in Fluids 68~(8) (2012) 1031--1052.

\bibitem{booy1978geometry}
M.~Booy, Geometry of fully wiped twin-screw equipment, Polymer Engineering \&
  Science 18~(12) (1978) 973--984.

\bibitem{sato2004heat}
S.~Sato, K.~Oka, A.~Murakami, Heat transfer behavior of melting polymers in
  laminar flow field, Polymer Engineering \& Science 44~(3) (2004) 423--432.

\bibitem{kalyon2007integrated}
D.~Kalyon, M.~Malik, An integrated approach for numerical analysis of coupled
  flow and heat transfer in co-rotating twin screw extruders, International
  Polymer Processing 22~(3) (2007) 293--302.

\bibitem{eusterholz2018cfd}
S.~Eusterholz, S.~Elgeti, Cfd-based optimization in plastics extrusion, in: AIP
  Conference Proceedings, Vol. 1960, AIP Publishing LLC, 2018, p. 110004.

\end{thebibliography}

\end{document}